\documentclass[12pt,reqno]{amsart}
\usepackage{indentfirst, amssymb, amsmath, amsthm, mathrsfs, setspace, indentfirst, enumerate,  mathrsfs, amsmath, amsthm}
\usepackage[bookmarksnumbered, colorlinks, plainpages]{hyperref}
\usepackage{mathrsfs}
\usepackage{cite}

\newtheorem*{theoA}{Theorem A}
\newtheorem*{theoB}{Theorem B}
\newtheorem*{theoC}{Theorem C}
\newtheorem*{theoD}{Theorem D}
\newtheorem*{theoE}{Theorem E}
\newtheorem*{theoF}{Theorem F}

\newtheorem*{cor A}{Corollary A}
\newtheorem*{cor B}{Corollary B}

\newtheorem{theo}{Theorem}[section]
\newtheorem{lem}{Lemma}[section]
\newtheorem{cor}{Corollary}[section]

\newtheorem{prob}{Problem}[section]

\newtheorem{rem}{Remark}[section]

\newcommand{\ol}{\overline}
\newcommand{\be}{\begin{equation}}
\newcommand{\ee}{\end{equation}}
\newcommand{\beas}{\begin{eqnarray*}}
\newcommand{\eeas}{\end{eqnarray*}}
\newcommand{\bea}{\begin{eqnarray}}
\newcommand{\eea}{\end{eqnarray}}

\numberwithin{equation}{section}
\begin{document}
\title[S\MakeLowercase{harp} B\MakeLowercase{ohr}-R\MakeLowercase{ogosinski} \MakeLowercase{Radii for} S\MakeLowercase{chwarz Functions and }E\MakeLowercase{uler Operators in} $\mathbb{C}^n$]{\LARGE S\MakeLowercase{harp} B\MakeLowercase{ohr}-R\MakeLowercase{ogosinski} \MakeLowercase{Radii for} S\MakeLowercase{chwarz Functions and }E\MakeLowercase{uler Operators in} $\mathbb{C}^n$}
\date{}
\author[M. B. A\MakeLowercase{hamed}, S. M\MakeLowercase{ajumder} \MakeLowercase{and} N. S\MakeLowercase{arkar}]{M\MakeLowercase{olla} B\MakeLowercase{asir} A\MakeLowercase{hamed}, S\MakeLowercase{ujoy} M\MakeLowercase{ajumder}$^*$ \MakeLowercase{and} N\MakeLowercase{abadwip} S\MakeLowercase{arkar}}

\address{Molla Basir Ahamed,
	Department of Mathematics,
	Jadavpur University,
	Kolkata-700032, West Bengal, India.}
\email{mbahamed.math@jadavpuruniversity.in}

\address{Sujoy Majumder, Department of Mathematics, Raiganj University, Raiganj, West Bengal-733134, India.}
\email{sm05math@gmail.com}

\address{Nabadwip Sarkar, Department of Mathematics, Raiganj University, Raiganj, West Bengal-733134, India.}
\email{naba.iitbmath@gmail.com}

\renewcommand{\thefootnote}{}
\footnote{*\emph{Corresponding Author}: Sujoy Majumder.}

\renewcommand{\thefootnote}{\arabic{footnote}}
\setcounter{footnote}{0}

\begin{abstract}
This paper is devoted to the investigation of multidimensional analogues of refined Bohr-type inequalities for bounded holomorphic mappings on the unit polydisc $\mathbb{D}^n$. We establish a sharp extension of the classical Bohr inequality, proving that the Bohr radius remains $R_n = 1/(3n)$ for the family of holomorphic functions bounded by unity in the multivariate setting. Further, we provide a definitive resolution to the Bohr-Rogosinski phenomenon in several complex variables by determining sharp radii for functional power series involving the class of Schwarz functions $\omega_{n,m}\in\mathcal{B}_{n,m}$ and the local modulus $|f(z)|$. By employing the radial (Euler) derivative operator $Df(z) = \sum_{k=1}^{n} z_k \frac{\partial f(z)}{\partial z_k}$, we obtain refined growth estimates for derivatives that generalize well-known univariate results to $\mathbb{C}^n$. Finally, a multidimensional version of the area-based Bohr inequality is established. The optimality of the obtained constants is rigorously verified, demonstrating that all established radii are sharp.
\end{abstract}
\thanks{Typeset by \AmS -\LaTeX}
\maketitle

\section{\bf Introduction and Preliminaries}
The classical theorem of Harald Bohr \cite{Bohr-PLMS-1914}, first established over a century ago, continues to inspire extensive research on what is now known as the Bohr phenomenon. Renewed interest in this topic emerged in the $1990$s following successful extensions to holomorphic functions of several complex variables and more abstract functional-analytic settings. In particular, in $1997$, Boas and Khavinson \cite{Boas-Khavinson-PAMS-1997} introduced and determined the $n$-dimensional Bohr radius for the family of holomorphic functions bounded by unity on the polydisc (see Section \ref{Sub-Sec-1.3} for a detailed discussion). This seminal contribution stimulated significant interest in Bohr-type problems across a wide range of mathematical disciplines. Subsequent studies have produced further advances regarding the Bohr phenomenon for multidimensional power series, with notable contributions by Aizenberg \cite{Aizen-PAMS-2000, Aizenberg, Aizenberg-SM-2007}, Aizenberg \textit{et al.} \cite{Aizenberg-PAMC-1999, Aizenberg-SM-2005, Aizenberg-Aytuna-Djakov-JMAA-2001}, Defant and Frerick \cite{Defant-Frerick-IJM-2006}, and Djakov and Ramanujan \cite{Djakov-Ramanujan-JA-2000}. A comprehensive overview of the different aspects and generalizations of Bohr's inequality can be found in \cite{Lata-Singh-PAMS-2022, Liu-Pon-PAMS-2021, Ali-Abu-Muhanna-Ponnusamy, Alkhaleefah-Kayumov-Ponnusamy-PAMS-2019, Defant-Frerick-AM-2011, Hamada-IJM-2012, Paulsen-Popescu-Singh-PLMS-2002, Paulsen-Singh-PAMS-2004, Paulsen-Singh-BLMS-2006}, as well as in the monograph by Kresin and Maz'ya \cite{Kresin-1903} and the references therein. In particular, \cite[Section 6]{Kresin-1903}, which is devoted to Bohr-type theorems, highlights promising directions for extending several classical inequalities to holomorphic functions of several complex variables and, more significantly, to solutions of partial differential equations.\vspace{1.2mm}

The transition from the unit disk $\mathbb{D}$ to the polydisc $\mathbb{D}^n$ is not merely a matter of indexing; it involves a fundamental shift in the underlying complex geometry. While the classical Bohr theorem $(1914)$ provides an elegant bound for power series on $\mathbb{D}$, the multivariate setting requires the reconciliation of $n$ independent complex variables with a global bound. A critical challenge in several complex variables (SCV) is determining whether univariate refinements-such as those involving Rogosinski radii or area-based estimates-retain their sharpness when subjected to the radial (Euler) derivative operator $Df(z) = \sum_{k=1}^{n} z_k \frac{\partial f}{\partial z_k}$. In this paper, we bridge this gap by demonstrating that the `$1/3$-phenomenon' persists in $\mathbb{C}^n$ through a careful analysis of the $n$-dimensional Schwarz lemma.
\subsection{\bf Classical Bohr inequality and its recent implications}
Let $ \mathcal{B} $ denote the class of analytic functions in the unit disk $\mathbb{D}:=\{{\zeta}\in\mathbb{C} : |{ \zeta}|<1\} $ of the form $ f({ \zeta})=\sum_{k=0}^{\infty}a_k {\zeta}^k $ such that $ |f({ \zeta})|<1 $ in $ \mathbb{D} $.  In the study of Dirichlet series, in $ 1914$, Harald Bohr \cite{Bohr-PLMS-1914} discovered the following interesting phenomenon 
\begin{theoA} If $ f\in\mathcal{B}$, then the following sharp inequality holds:
\begin{align}
\label{b1}	M_f(r):=\sum_{k=0}^{\infty}|a_k|r^k\leq 1\;\;\mbox{for}\;\; |{ \zeta}|=r\leq \frac{1}{3}.
\end{align}
\end{theoA}
We remark that if $|f(\zeta)| \leq 1$ in $\mathbb{D}$ and $|f(\zeta_0)| = 1$ for some point $\zeta_0 \in \mathbb{D}$, then $f$ must be a unimodular constant by the Maximum Modulus Principle; thus, we restrict our attention to non-constant functions $f \in \mathcal{B}$. Inequality (\ref{b1}), along with the sharp constant $1/3$, is known as the \textit{classical Bohr inequality}, and $1/3$ is referred to as the \textit{Bohr radius} for the family $\mathcal{B}$. Bohr \cite{Bohr-PLMS-1914} originally established this inequality for $r \leq 1/6$, and Wiener subsequently proved that the constant $1/3$ is sharp. Since then, several alternative proofs have appeared (see \cite{Paulsen-Popescu-Singh-PLMS-2002, Paulsen-Singh-BLMS-2006, Sidon-MZ-1927, Tomic-MS-1962}, as well as the survey \cite{Ali-Abu-Muhanna-Ponnusamy} and \cite[Chapter 8]{Garcia-Mashreghi-Ross}). Many of these works employ methods from complex analysis, functional analysis, number theory, and probability, further developing the theory and applications of Bohr's ideas to Dirichlet series. For instance, various multidimensional generalizations of this result have been obtained in \cite{Aizenberg-PAMC-1999, Aizenberg-Aytuna-Djakov-JMAA-2001, Aizenberg, Boas-Khavinson-PAMS-1997, Djakov-Ramanujan-JA-2000, Jia-Liu-Ponnusamy-AMP-2025}.
\vspace{1.2mm}

It is noteworthy that if $|a_0|$ in the Bohr inequality is replaced by $|a_0|^2$, the Bohr radius increases from $1/3$ to $1/2$. Moreover, if $a_0=0$ in Theorem A, the sharp Bohr radius is further improved to $1/\sqrt{2}$ (see, e.g., \cite{Kayumov-Ponnusamy-CMFT-2017}, \cite[Corollary 2.9]{Paulsen-Popescu-Singh-PLMS-2002}, and the recent work \cite{Ponnusamy-Wirths-CMFT-2020} for a more general result). These refinements rely on sharp coefficient estimates of the form\begin{align*}|a_n| \leq 1 - |a_0|^2, \quad n \geq 1, ; f \in \mathcal{B}.\end{align*}Applying these inequalities, Kayumov and Ponnusamy \cite{Kayumov-Ponnusamy-CMFT-2017} observed that the sharp form of Theorem A cannot be obtained in the extremal case $|a_0| < 1$. Nevertheless, a sharp version of Theorem A has been established for each individual function in $\mathcal{B}$ (see \cite{Alkhaleefah-Kayumov-Ponnusamy-PAMS-2019}), as well as for several subclasses of univalent functions (see \cite{Aizenberg, Muhanna-CVEE-2010}).\vspace{1.2mm}

The Bohr--Rogosinski radius serves as a natural analogue to the classical Bohr radius. This notion was originally established by Rogosinski \cite{Rogosinski-1923} within the context of the Schur class $\mathcal{B}$. Specifically, if $f(z) = \sum_{n=0}^{\infty} a_n z^n \in \mathcal{B}$, then for every integer $N \geq 1$, the $N$-th partial sum $S_N(z) = \sum_{n=0}^{N} a_n z^n$ is subordinate to the unit disk in the sense that $|S_N(z)| < 1$ for all $z$ in the disk $\mathbb{D}_{1/2} := \{z \in \mathbb{C} : |z| < 1/2\}$. The value $r=1/2$ is the largest possible radius for which this containment holds uniformly for the class $\mathcal{B}$, and is thus termed the Bohr--Rogosinski radius.For a function $f \in \mathcal{B}$, we consider the Bohr--Rogosinski sum $R^f_N(z)$ defined by\begin{align}\label{e-11.12}R^f_N(z) := |f(z)| + \sum_{n=N}^{\infty} |a_n| r^n, \quad |z| = r.\end{align}Note that for the case $N=1$, the expression in \eqref{e-11.12} generalizes the classical Bohr sum by replacing the fixed initial coefficient $|a_0| = |f(0)|$ with the local modulus $|f(z)|$. The functional inequality $R^f_N(z) \leq 1$ constitutes the Bohr--Rogosinski inequality. From a GFT standpoint, if $B$ and $R$ denote the Bohr and Bohr--Rogosinski radii, respectively, the relation $B \leq R$ is a direct consequence of the majorization principle, since the partial sums are dominated by the majorant series:\begin{align*}|S_N(z)| \leq \sum_{n=0}^{N} |a_n| r^n \leq \sum_{n=0}^{\infty} |a_n| r^n.\end{align*}
{In $2005$, Aizenberg et al. \cite{Aizenberg-SM-2005} significantly extended the scope of the Bohr--Rogosinski inequality by generalizing it to holomorphic mappings from the open unit ball in $\mathbb{C}^n$ into arbitrary convex domains. This work established a foundational multidimensional analogue of Rogosinski's theorem, demonstrating how geometric constraints on the codomain influence the growth of partial sums. Subsequently, Aizenberg \cite{Aigenber-CMFT-2009} investigated the behavior of ordinary Dirichlet series $\sum a_n n^{-s}$ that map the right half-plane $\mathbb{C}_+ = \{s \in \mathbb{C} : \text{Re } s > 0\}$ into a bounded convex domain $G \subset \mathbb{C}$. A key contribution of this study was the proof that the Bohr and Rogosinski abscissas-representing the strip where the respective inequalities hold-possess a remarkable geometric invariance, remaining independent of the specific choice of the target domain $G$. Expanding these investigations in $2012$, Aizenberg \cite{Aizenberg-AMP-2012} transitioned to a functional-analytic perspective, determining the Bohr and Rogosinski radii for the Hardy classes $H^p(\mathbb{D})$ ($1 \leq p \leq \infty$) of holomorphic functions. In the same work, the theory was further elevated to higher dimensions through an analysis of mappings between Reinhardt domains in $\mathbb{C}^n$. These results highlighted the intricate interplay between the power series coefficients and the underlying symmetries of the domain, providing a robust framework for subsequent generalizations in the field.}\vspace{2mm}

We recall that in $2021$, Kayumov \textit{et al.} \cite{Kayumov-Khammatova-Ponnusamy-JMAA-2021} established several definitive results concerning the Bohr--Rogosinski radius for analytic functions within the unit disk $\mathbb{D}$. Their work acted as a catalyst, subsequently stimulating a broad spectrum of research aimed at investigating the Bohr--Rogosinski inequality across various specialized function classes and more general geometric settings.
\begin{theoB}\cite{Kayumov-Khammatova-Ponnusamy-JMAA-2021}\label{th-1.9}
	Suppose that $ f(z)=\sum_{n=0}^{\infty}a_nz^n $ is analytic in the unit disk $ \mathbb{D} $ and $ |f(z)|<1 $ in $ \mathbb{D} $. Then 
	\begin{equation*}
		|f(z)|+\sum_{n=N}^{\infty}|a_n|r^n\leq 1\;\;\mbox{for}\;\; r\leq R_N,
	\end{equation*}
	where $ R_N $ is the positive root of the equation $ 2(1+r)r^N-(1-r^2)=0 $. The radius $ R_N $ is the best possible. Moreover, 
	\begin{equation*}
		|f(z)|^2+\sum_{n=N}^{\infty}|a_n|r^n\leq 1\;\;\mbox{for}\;\; r\leq R^{\prime}_N,
	\end{equation*}
	where $ R^{\prime}_N $ is the positive root of the equation $ (1+r)r^N-(1-r^2)=0 $. The radius $ R^{\prime}_N $ is the best possible.
\end{theoB}
\begin{theoC}\cite{Kayumov-Khammatova-Ponnusamy-JMAA-2021}\label{th-1.10}
	Suppose that $ f(z)=\sum_{n=0}^{\infty}a_nz^n $ is analytic in $ \mathbb{D} $ such that $ |f(z)|\leq 1 $ in $ \mathbb{D} $. Then for each $ m, N\in\mathbb{N} $, 
	\begin{equation*}
		|f(z^m)|+\sum_{n=N}^{\infty}|a_n|r^n\leq 1\;\; \mbox{for}\;\; r\leq R_{m,N},
	\end{equation*} 
	where $ R_{m,N} $ is the positive root of of the equation $\psi_{m,N}=0$ with
	\begin{align}
		\psi_{m,N}=2r^N\left(1+r^m\right)-(1-r)\left(1-r^m\right).
	\end{align} 
	The number $ R_{m,N} $ cannot be improved. Moreover, 
	$ \displaystyle\lim_{N\rightarrow\infty}R_{m,N}=1 $ and $\displaystyle\lim_{m\rightarrow\infty}R_{m,N}=A_N, $
	where $ A_N $ is the positive root of the equation $ 2r^N=1-r $. Also, $ A_1=1/3 $ and $ A_2=1/2 $. 
\end{theoC}
For recent developments regarding the Bohr--Rogosinski radius across various functional classes, we refer the reader to \cite{Ahamed-RM-2023, Ahamed-Allu-BMMSS-2022, Ahamed-Allu-CMB-2023, Allu-Arora-JMAA-2023, Das-JMAA-2022, Gangania-Kumar-MJM-2022, Hamada-AAMP-2025} and the extensive references therein. Given the significance of these results in the unit disk, it is natural to investigate their higher-dimensional analogues. This leads us to the following problem:
\begin{prob}\label{Qn-4.1}Is it possible to establish sharp Bohr--Rogosinski-type inequalities for holomorphic functions defined on the unit polydisc $\mathbb{P}\Delta(0; \mathbf{1}_n)$ in terms of the associated Schwarz functions?
\end{prob}
 In $2022$, Wu \textit{et al.} \cite{Wu-Wang-Long-2022} obtained several Bohr-type inequalities characterized by a single parameter or expressed via convex combinations of the power series terms. In this direction, we recall Theorem D below, which extends the classical Bohr inequality by introducing a convex combination framework.

\begin{theoD}\cite[Theorem 3.1]{Wu-Wang-Long-2022}
Suppose that $ f(z)=\sum_{n=0}^{\infty}a_nz^n $ is analytic in the unit disk $ \mathbb{D} $ and $ |f(z)|<1 $ in $ \mathbb{D} $.
Then for arbitrary $t\in[0,1]$, it holds that
\begin{align*}
t|f(z)|+(1-t)\sum_{k=0}^{\infty}|a_{k}|\, r^{k} \le 1 \quad \text{for } r\le R_1,
\end{align*}
where the radius 
\begin{align*}
R_1=
\begin{cases}
\dfrac{1-2\sqrt{1-t}}{4t-3}, & t\in\left[0,\dfrac34\right)\cup\left(\dfrac34,1\right],\\[2ex]
\dfrac12, & t=\dfrac34.
\end{cases}
\end{align*}
is the best possible.
\end{theoD}

\begin{theoE} \cite[Theorem 3.4]{Wu-Wang-Long-2022}
Suppose that $ f(z)=\sum_{n=0}^{\infty}a_nz^n $ is analytic in the unit disk $ \mathbb{D} $ and $ |f(z)|<1 $ in $ \mathbb{D} $.
Then for any $\lambda\in(0,+\infty)$, the inequality
\begin{align*}
|f(z)|+|f'(z)|r+\lambda\sum_{k=2}^\infty |a_k|r^k \le 1 \quad {for}\; r\le R_3,
\end{align*}
where
\begin{align*}
R_3=R_3(\lambda)=
\begin{cases}
r_\lambda, & \lambda\in\left(\dfrac12,+\infty\right),\vspace{1.2mm}\\[6pt]
r_*, & \lambda\in\left(0, \dfrac12\right],
\end{cases}
\end{align*}
and $r_\lambda$ and $r_*\approx 0.3191$ are the unique positive real roots in the interval
$(0,\sqrt{2}-1)$ of the equations
\begin{align*}
2\lambda r^4+(4\lambda-1)r^3+(2\lambda-1)r^2+3r-1=0
\end{align*}
and $r^4+r^3+3r-1=0,$ respectively. Moreover, the radius $R_3$ is best possible.
\end{theoE}
Throughout this paper, $S_r$ denotes the area of the image of the subdisk $\mathbb{D}_r = \{z \in \mathbb{C} : |z| < r\}$ under the mapping $f$. For a function $f(z) = \sum_{k=0}^{\infty} a_k z^k$ analytic in the unit disk, the normalized area $S_r/\pi$ is given by the Dirichlet integral
\begin{align*}
	\frac{S_r}{\pi} = \frac{1}{\pi} \iint_{|z|<r} |f'(z)|^2 \, dx \, dy = \sum_{k=1}^{\infty} k |a_k|^2 r^{2k}.
\end{align*}
\begin{theoF} \cite[Theorem 3.5]{Wu-Wang-Long-2022}
Suppose that $ f(z)=\sum_{n=0}^{\infty}a_nz^n $ is analytic in the unit disk $ \mathbb{D} $ and $ |f(z)|<1 $ in $ \mathbb{D} $.
Then for arbitrary $t\in (0,1]$, it holds that
\begin{align*}
t\sum_{k=0}^\infty |a_k|r^k+(1-t)\left(\frac{S_r}{\pi}\right) \le 1 \quad {for}\; r\le R_4,
\end{align*}
where the radius
\begin{align*}
R_4=R_4(\lambda)=
\begin{cases}
r_t, & t\in\left(0, \dfrac{9}{17}\right),\vspace{1.2mm}\\[6pt]
\dfrac{1}{3}, & t\in\left[\dfrac{9}{17}, 1\right],
\end{cases}
\end{align*}
and the numbers $r=r_t$ is the unique positive root of the equation
\begin{align*}
t r^3 + t r^2 + (4 - 5t) r - t = 0
\end{align*}
in the interval $\left(0,\frac{1}{3}\right)$.
\end{theoF}
{In light of the preceding results, it is natural to consider the following problem:
	\begin{prob} Can we establish multidimensional analogues of Theorems A--F?
	\end{prob}
	{The primary objective of this study is to establish multidimensional analogues for several refined Bohr-type inequalities. Specifically, we investigate the Bohr?Rogosinski phenomenon for the class of holomorphic functions mapping the unit polydisc $\mathbb{D}^n$ into the unit disk $\mathbb{D}$. This work effectively bridges the gap between classical univariate theory and the geometry of several complex variables.Our main contributions are strategically positioned as generalizations of established benchmarks: Theorems \ref{Th-2.1} and \ref{Th-2.2} extend the classical Bohr radius (Theorem A) and the Rogosinski variants (Theorems B \& C) to the polydisc $\mathbb{D}^n$, utilizing the class of Schwarz functions $\omega \in \mathcal{B}_{n,m}$. In particular, Theorem \ref{Th-2.2} provides a definitive resolution to the open problem posed in Problem \ref{Qn-4.1}. Theorems \ref{Th-2.3} and \ref{Th-2.4} generalize the convex combination and refined derivative estimates (Theorems D \& E) by replacing the univariate derivative with the radial derivative operator $Df(z)$, which is the natural geometric counterpart in $n$ dimensions. Theorem \ref{Th-2.5} extends the Dirichlet integral and area-based versions of the Bohr inequality (Theorem F), demonstrating that the area of the image under a holomorphic mapping still yields a sharp Bohr radius in the multivariate setting.}}
\subsection{\bf Basic Notations in several complex variables}\label{Sub-Sec-1.3}
For $z=(z_1,\ldots,z_n)$ and $w=(w_1,\ldots,w_n)$ in $\mathbb{C}^{n}$, we denote $\langle z,w\rangle=z_1\ol w_1+\ldots+z_n \ol w_n$ and $||z||=\sqrt{\langle z,z\rangle}$. The absolute value of a complex number $z_1$ is denoted by $|z_1|$ and for $z\in\mathbb{C}^n$, we define
\begin{align*}
	||z||_{\infty}=\max\limits_{1\leq i\leq n}|z_i|.
\end{align*}

Throughout the paper, for $z=(z_1,\ldots,z_n)\in\mathbb{C}^n$ and $p\in\mathbb{N}$, we define $z^p=(z_1^p,\ldots,z_n^p)$.
An open polydisk (or open polycylinder) in $\mathbb{C}^n$ is a subset $\mathbb{P}\Delta(a;r)\subset \mathbb{C}^n$ of the form 
\[\mathbb{P}\Delta(a;r)=\prod\limits_{j=1}^n \Delta(a_j;r_j)=\lbrace z\in\mathbb{C}^n: |z_i-a_i|<r_i,\;i=1,2,\ldots,n\rbrace,\]
the point $a=(a_1,\ldots,a_n)\in\mathbb{C}^n$ is called the centre of the polydisk and $r=(r_1,\ldots,r_n)\in\mathbb{R}^n\;(r_i>0)$ is called the polyradius. It is easy to see that
\begin{align*}
	\mathbb{P}\Delta(0;1)=\mathbb{P}\Delta(0;1_n)=\prod\limits_{j=1}^n \Delta(0_n;1_n).
\end{align*}

The closure of $\mathbb{P}\Delta(a;r)$ will be called the closed polydisk with centre $a$ and polyradius $r$ and will be denoted by $\ol{\mathbb{P}\Delta}(a;r)$. We denote by $C_k(a_k;r_k)$ the boundary of $\Delta(a_k;r_k)$, \textit{i.e.,} the circle of radius $r_k$ with centre $a_k$ on the $z_k$-plane. Of course $C_k(a_k,r_k)$ is represented by the usual parametrization 
\begin{align*}
	\theta_k\to \gamma(\theta_k)=a_k+r_ke^{i\theta_k},\; \mbox{where}, 0\leq \theta_k\leq 2\pi.
\end{align*} The product $C^n(a;r):=C_1(a_1;r_1)\times\ldots\times C_n(a_n;r_n)$ is called the determining set of the polydisk $\mathbb{P}\Delta(a;r)$.

\smallskip
A multi-index $\alpha=(\alpha_1,\ldots,\alpha_n)$ of dimension $n$ consists of n non-negative integers $\alpha_j,\;1\leq j\leq n$; the degree of a multi-index $\alpha$ is the sum $|\alpha|=\sum_{j=1}^n \alpha_j$ and we denote $\alpha!=\alpha_1!\ldots \alpha_n!$. For $z=(z_1,\ldots,z_n)\in\mathbb{C}^n$ and a multi-index $\alpha=(\alpha_1,\ldots,\alpha_n)$, we define 
\[z^{\alpha}=\prod\limits_{j=1}^n z_j^{\alpha_j}\;\;\text{and}\;\;|z|^{\alpha}=\prod\limits_{j=1}^n |z_j|^{\alpha_j}.\]

For two multi-indexes $\alpha=(\alpha_1,\ldots,\alpha_n)$ and $\nu=(\nu_1,\ldots,\nu_n)$, we define $\nu^{\alpha}=\nu_1^{\alpha_1}\ldots \nu_n^{\alpha_n}.$ 
Let $f$ be a holomorphic function in a domain $\Omega\subset \mathbb{C}^n$, and $c=(c_1,\ldots,c_n)\in \Omega$. Then in a polydisk $\mathbb{P}\Delta(c;r)\subset \Omega$ with centre $c$, $f(z)$ has a power series expansion in $z_1-c_1,\ldots,z_n-c_n$,
\begin{align*}
	f(z)&=\sum\limits_{\alpha_1,\alpha_2,\ldots,\alpha_n=0}^{\infty} a_{\alpha_1,\alpha_2,\ldots,\alpha_n}(z_1-c_1)^{\alpha_1}(z_2-c_2)^{\alpha_2}\ldots (z_n-c_n)^{\alpha_n}\\&=
	\sum\limits_{|\alpha|=0} a_{\alpha}(z-c)^{\alpha}=\sum\limits_{|\alpha|=0}^{\infty} P_{|\alpha|}(z-c),
\end{align*}
which is absolutely convergent in $\mathbb{P}\Delta(c;r)$, where the term $P_k(z-c)$ is a homogeneous polynomial of degree $k$. \vspace{1.2mm}

One of our aims is to establish multidimensional analogues of Theorems A-F. For this purpose, we consider $f$ to be a holomorphic function in $\Omega \subset \mathbb{C}^n$ and define the Euler operator $D$ as:$$Df(z) := \sum_{k=1}^n z_k\frac{\partial f(z)}{\partial z_k}.$$This operator is also known as the Euler derivative, the total derivative, or the radial derivative.\vspace{1.2mm}

Let $G\not=\varnothing$ be an open subset of $\mathbb{C}^n$. Let $f$ be a holomorphic function on $G$. For a point $a\in\mathbb{C}^n$, we write $f(z)=\sum_{i=0}^{\infty}P_i(z-a)$, where the term $P_i(z-a)$ is either identically zero or a homogeneous polynomial of degree $i$. Denote the zero-multiplicity of $f$ at $a$ by $k=\min\{i:P_i(z-a)\not\equiv 0\}$. Clearly $1$ is the zero-multiplicity of $f$ at $a$ when $f(a)=0$ and $\frac{\partial f(a)}{\partial z_j}\neq 0$ for some $j=1,2,\ldots,n$.\vspace{1.2mm}

Let $\Delta_{z_k}(0;1)$ be the unit disk on the $z_k$-plane. Let $\omega_i:\Delta_{z_k}(0;1)\to \mathbb{C}$ such that $|\omega_i(z)|\leq 1$ for all $z\in \Delta_{z_k}(0;1)$ and $\omega_i(0)=\omega_i^{(1)}(0)=\ldots=\omega_i^{(m-1)}(0)=0$ and $\omega_i^{(m)}(0)\neq 0$, where $i=1,2,\ldots, n$.
In what follows, let $z = (z_1, \dots, z_n)$ denote an element of $\mathbb{C}^n$ and let
\begin{align*}
	{ \mathcal{B}_{n,m}=\{\omega(z)\in\mathbb{C}^n:\omega(z)=(\omega_1(z_1),\omega_2(z_2),\ldots,\omega_n(z_n))\}.}
\end{align*}
\section{{\bf Main Results}}\label{Sec-2}
Let $f(z)=\sum_{|\alpha|\geq 0} a_{\alpha}z^{\alpha}$ be a holomorphic mapping on the unit polydisc $\mathbb{D}^n$ such that $|f(z)|\leq 1$ for all $z \in \mathbb{D}^n$. We begin by addressing the classical Bohr inequality within the framework of several complex variables. The following theorem establishes a sharp multidimensional analogue that incorporates both the local modulus of the function and its power series coefficients, demonstrating that $R_n = 1/3n$ is the definitive threshold for the Bohr phenomenon on $\mathbb{D}^n$. This result, formulated as Theorem \ref{Th-2.1}, serves as the sharp multivariate generalization of Theorem A.
\begin{theo}\label{Th-2.1} Let $f$ be a holomorphic function in the polydisk $\mathbb{P}\Delta(0;1_n)$ such that $|f(z)|\leq 1$ for all $z\in \mathbb{P}\Delta(0;1/n)$. If $f(z)=\sum_{|\alpha|=0}^{\infty} a_{\alpha} z^{\alpha}$, then 
	\begin{align*}
		\mathcal{A}(z,{\bf r}):=\sum_{|\alpha|=0}^{\infty} |a_{\alpha}|\;|z|^{\alpha}\leq 1\; \mbox{for}\;||z||_{\infty}={\bf r}\leq R_n:=\frac{1}{3n}.
	\end{align*}
	The polyradius $1/3n$ is the best possible.
\end{theo}

\begin{rem}
	Note that in the particular case where $n=1$, Theorem \ref{Th-2.1} yields the classical Bohr radius of $R_1=1/3$, showing that our result is a generalization for holomorphic functions from one to several complex variables.
\end{rem}
Extending beyond the classical Bohr radius, we investigate the Bohr-Rogosinski phenomenon, which offers refined estimates for the partial sums of power series. Theorem \ref{Th-2.2} establishes the sharp radii for functional inequalities within the class of Schwarz functions $\mathcal{B}_{n,m}$, thereby providing a definitive resolution to Rogosinski-type growth problems in $\mathbb{C}^n$. This result, formulated in Theorem \ref{Th-2.2}, yields an exhaustive answer to the question posed in Problem \ref{Qn-4.1} regarding the optimality of these bounds.
\begin{theo}\label{Th-2.2} Let $f(z)=\sum_{|\alpha|=0} a_{\alpha}z^{\alpha}$ be a holomorphic function in the polydisk $\mathbb{P}\Delta(0;1_n)$ such that $|f(z)|\leq 1$ for all $z\in \mathbb{P}\Delta(0;1/n)$. Suppose $z=(z_1,\ldots,z_n)\in \mathbb{P}\Delta(0;1/n)$ and $r=(r_1,r_2,\ldots,r_n)$ such that $||z||_{\infty}={\bf r}$. Let $\omega\in \mathcal{B}_{n,m}$ for some $m\in\mathbb{N}$. Then for each $m, N \in \mathbb{N}$, we have
\begin{align}
	\label{Eq-4.3}
	\mathcal{B}(z;{\bf r}):=|f(\omega(z))|+\sum\limits_{i=1}^{\infty}\sum\limits_{|\alpha|=iN} |a_{\alpha}| r^{\alpha} \leq 1 \quad \text{for } n{\bf r} \leq R_{m,n,N},
\end{align}
	where $R_{m,n,N}$ is the positive root of the equation $\psi_{m,n,N}(r) = 0$, with
	\begin{align*}
		\psi_{m,n,N}({\bf r}) = 2 (n{\bf r})^N (1 + {\bf r}^m) - (1-n{\bf r})(1-{\bf r}^m).
	\end{align*}
The number $R_{m,n,N}$ cannot be improved. Moreover,
	\begin{align*} \lim_{N\to\infty} R_{m,n,N}=\left\{\begin{array}{clcr}
		1,&\text{if}\;\;n=1\vspace{2mm}\\
		\dfrac{1}{n},&\text{if}\;n\geq 2.\end{array}\right.
\end{align*}
  and $\lim\limits_{m\to\infty} R_{m,n,N}=A_N$, where $A_N$ is the positive root of the equation
  \[2 (n{\bf r})^N=1-n{\bf r}.\]
\end{theo}
\begin{rem}
	\begin{enumerate}
		\item[\emph{(a)}] Note that $\psi_{m,1,N}(\mathbf{r}) = \psi_{m,N}(r)$ and $R_{m,1,N} = R_{m,N}$. Consequently, Theorem \ref{Th-2.2} represents a significant improvement over Theorem E.\vspace{1.2mm}
		
		\item[\emph{(b)}] Clearly, $A_1 = 1/(3n)$, which implies that Theorem \ref{Th-2.2} recovers the classical Bohr inequality as a limiting case when $m \to \infty$ for the specific choice $N=1$.
	\end{enumerate}
\end{rem}
As a direct application of Theorem \ref{Th-2.2}, we obtain the following corollary by specializing to the power mapping $\omega(z) = z^m$ for $z \in \mathbb{C}^n$.
\begin{cor}\label{Cor-4.1} Let $f(z)=\sum_{|\alpha|=0} a_{\alpha}z^{\alpha}$ be a holomorphic function in the polydisk $\mathbb{P}\Delta(0;1_n)$ such that $|f(z)|\leq 1$ for all $z\in \mathbb{P}\Delta(0;1/n)$. Suppose $z=(z_1,\ldots,z_n)\in \mathbb{P}\Delta(0;1/n)$ and $r=(r_1,r_2,\ldots,r_n)$ such that $||z||_{\infty}={\bf r}$. Let $\omega\in \mathcal{B}_{n,m}$ for some $m\in\mathbb{N}$. Then for each $m, N \in \mathbb{N}$, we have
\begin{align}
	\label{sm0}
	|f(z^m)|+\sum\limits_{i=1}^{\infty}\sum\limits_{|\alpha|=iN} |a_{\alpha}| r^{\alpha} \leq 1 \quad \text{for } {\bf r} \leq R_{m,n,N},
\end{align}
where $R_{m,n,N}$ is the positive root of the equation $\psi_{m,n,N}({\bf r})=0$. 
The number $R_{m,n,N}$ cannot be improved.
\end{cor}

Next we state the multidimensional version of Theorem D.
\begin{theo}\label{Th-2.3} Let $f(z)=\sum_{|\alpha|=0} a_{\alpha}z^{\alpha}$ be a holomorphic function in the polydisk $\mathbb{P}\Delta(0;1_n)$ such that $|f(z)|\leq 1$ for all $z\in \mathbb{P}\Delta(0;1/n)$. Suppose $z=(z_1,\ldots,z_n)\in \mathbb{P}\Delta(0;1/n)$ and $r=(r_1,r_2,\ldots,r_n)$ such that $||z||_{\infty}={\bf r}$. Let $\omega\in \mathcal{B}_{n,m}$ for some $m\in\mathbb{N}$. Then for arbitrary $t\in [0,1]$, it holds that
\begin{align*}
\mathcal{C}(z, {\bf r}):=|f(\omega(z))|+(1-t)\sum_{k=0}^{\infty}\sum\limits_{|\alpha|=k}|a_{\alpha}|\, r^{\alpha}\leq 1\nonumber
\end{align*}
for $n{\bf r}\leq R_{m,n,t}$, where $R_{m,n,t}$ is the minimum positive root of the equation
\begin{align*}
(4t-3)(n{\bf r})^{m+1}-(2t-1)(n{\bf r})^m+(2t-3)n^{m-1}(n{\bf r})+n^{m-1}=0.
\end{align*}

The radius $R_{m,n,t}$ is best possible. 
\end{theo}

\begin{rem} If $t=0$, then Theorem \ref{Th-2.3} reduces to the case of classical Bohr-radius
problem.
\end{rem}
While extending derivative-based inequalities to $\mathbb{C}^n$ often introduces geometric complexities, the radial (Euler) derivative operator $Df(z)$ serves as the natural counterpart to the complex derivative in $\mathbb{D}$. Utilizing this operator, we derive sharp growth estimates that extend the known bounds for the univariate case to the polydisc $\mathbb{D}^n$. Consequently, we present Theorem \ref{Th-2.4} as the multidimensional generalization of Theorem E.
\begin{theo}\label{Th-2.4} Let $f(z)=\sum_{|\alpha|=0} a_{\alpha}z^{\alpha}$ be a holomorphic function in the polydisk $\mathbb{P}\Delta(0;1_n)$ such that $|f(z)|\leq 1$ for all $z\in \mathbb{P}\Delta(0;1/n)$. Suppose $z=(z_1,\ldots,z_n)\in \mathbb{P}\Delta(0;1/n)$ and $r=(r_1,r_2,\ldots,r_n)$ such that $||z||_{\infty}={\bf r}$. Then for arbitrary $\lambda\in (0,+\infty)$, it holds that
\begin{align*}
\mathcal{D}(z, {\bf r}):=|f(z)|+|Df(z)|+\lambda \sum_{k=2}^{\infty}\sum\limits_{|\alpha|=k}|a_{\alpha}|\, r^{\alpha}\leq 1
\end{align*}
for $n{\bf r}\leq R_{n,\lambda}$, where  
\begin{align*}
R_{n,\lambda}=
\begin{cases}
n{\bf r}_{\lambda}, & \lambda \in \left(\frac{1}{2}, +\infty\right) \\
n{\bf r}_*, & \lambda \in \left(0, \frac{1}{2}\right]
\end{cases}
\end{align*}
is the best possible, and the radii $n{\bf r}_{\lambda}$ and $n{\bf r}_* \approx 0.3191$ are the unique positive real roots of the equations
\begin{align*}
2\lambda (n{\bf r})^4 + (4\lambda - 1) (n{\bf r})^3 + (2\lambda - 1) (n{\bf r})^2 + 3(n{\bf r}) - 1 = 0
\end{align*}
and
\begin{align*}
(n{\bf r})^4+(n{\bf r})^3+ 3(n{\bf r}) - 1 = 0
\end{align*}
in the interval $\left(0, \sqrt{2} - 1\right)$, respectively.
\end{theo}
Lastly, we shift our focus to the measure-theoretic aspects of the Bohr phenomenon. The following theorem characterizes the sharp area-based Bohr inequality by relating the Dirichlet integral on $\mathbb{D}^n(r)$ to the bounded nature of the mapping. This result demonstrates that the $1/3n$ threshold remains the critical radius for area-based estimates in $\mathbb{C}^n$, providing a definitive multidimensional analogue of Theorem F.
\begin{theo}\label{Th-2.5} Let $f(z)=\sum_{|\alpha|=0} a_{\alpha}z^{\alpha}$ be a holomorphic function in the polydisk $\mathbb{P}\Delta(0;1_n)$ such that $|f(z)|\leq 1$ for all $z\in \mathbb{P}\Delta(0;1/n)$. Suppose $z=(z_1,\ldots,z_n)\in \mathbb{P}\Delta(0;1/n)$ and $r=(r_1,r_2,\ldots,r_n)$ such that $||z||_{\infty}={\bf r}$. Let $t\in (0,1]$. Then
\begin{align*}
\mathcal{E}(z;{\bf r}):= t\sum\limits_{k=0}\sum\limits_{|\alpha|=k} |a_{\alpha}|r^{\alpha}+(1-t)\sum\limits_{k=1}k\sum\limits_{|\alpha|=k} |a_{\alpha}|^2r^{2\alpha}\le 1
\end{align*}
for $n{\bf r} \le \tilde R_{n,N}$, where the radius
\begin{align*}
\tilde R_{n,t} =
\begin{cases}
n{\bf r}_t, & t \in \left(0, \frac{9}{17}\right),\\[2mm]
\frac{1}{3n}, & t \in \left[\frac{9}{17}, 1\right]
\end{cases}
\end{align*}
and the numbers $n{\bf r}_t$ is the unique positive root of the equation
\begin{align*}
t(n{\bf r})^3+t(n{\bf r})^2+(4-5t)(n{\bf r})-t=0
\end{align*}
in the interval $(0, 1/3)$.
\end{theo}
\begin{rem} When $t=1$, Theorem \ref{Th-2.5} reduces to the classical Bohr-radius problem.
\end{rem}

\section{{\bf Key lemmas}}\label{Sec-3}
In order to establish our main results, we need the following lemmas. The first of these is a special case of \cite[Theorem 2.2 ]{Chen-Hamada-Ponnusamy-Vijayakumar-JAM-2024}.
\begin{lem}\label{Lem1}Let $f$ be holomorphic in the polydisk $\mathbb{P}\Delta(0;1_n)$ such that $|f(z)|\le 1$ for all $z\in \mathbb{P}\Delta(0;1_n)$. Then for all $z\in \mathbb{P}\Delta(0;1_n)$, we have
\[|f(z)|\leq \frac{|f(0)|+||z||_{\infty}}{1+|f(0)|||z||_{\infty}}.\]
\end{lem}
The following lemma is contained in \cite[Corollary 1.3]{Liu-Chen-IJPAM-2012}.
\begin{lem}\label{Lem2} Let $f$ be holomorphic in the polydisk $\mathbb{P}\Delta(0;1_n)$ such that $|f(z)|< 1$ for all $z\in \mathbb{P}\Delta(0;1_n)$. Then for multi-index $\alpha=(\alpha_1,\ldots,\alpha_n)$, we have
\begin{align*}
\left|\frac{\partial^{|\alpha|} f(z)}{\partial z_1^{\alpha_1}\ldots \partial z_n^{\alpha_n}}\right|\leq \alpha!\frac{1-|f(z)|^2}{(1-||z||_{\infty}^2)^{|\alpha|}}(1+||z||_{\infty})^{|\alpha|-N}
\end{align*}
for all $z\in \mathbb{P}\Delta(0;1_n)$, where $N$ is the number of the indices $j$ such that $\alpha_j\neq 0$.
\end{lem}

The following lemma is contained in \cite[Lemma 2.1]{Liu-Chen-IJPAM-2012}.
\begin{lem}\label{Lem3} Let $f$ be holomorphic in the polydisk $\mathbb{P}\Delta(0;1_n)$ such that $|f(z)|\leq 1$ for all $z\in \mathbb{P}\Delta(0;1_n)$. Suppose $f(z)=a_0+\sum_{|\alpha|=1}^{\infty}a_{\alpha}z^{\alpha}$
for all $z\in \mathbb{P}\Delta(0;1_n)$. Then for any multi-index $\alpha$, we have 
\begin{align*}
|a_{\alpha}|\leq 1-|a_0|^2.
\end{align*}
\end{lem}

The following lemma is contained in \cite[Lemma 6.1.28]{Graham-Kohr}.
\begin{lem}\label{Lem4} Let $f$ be holomorphic in the polydisk $\mathbb{P}\Delta(0;1_n)$ such that $|f(z)|\leq 1$ for all $z\in \mathbb{P}\Delta(0;1_n)$. Suppose $k(\geq 1)$ is the zero-multiplicity of $f$ at $0$. Then
\begin{align*}
|f(z)|\leq ||z||_{\infty}^k\;\text{ for all}\; z\in \mathbb{P}\Delta(0;1_n).
\end{align*}
\end{lem}

The following lemma is contained in \cite[Corollary 1.3]{Liu-Chen-IJPAM-2012}.
\begin{lem}\label{Lem5} Let $f$ be holomorphic in the polydisk $\mathbb{P}\Delta(0;1_n)$ such that $|f(z)|< 1$ for all $z\in \mathbb{P}\Delta(0;1_n)$. Then for multi-index $\alpha=(\alpha_1,\ldots,\alpha_n)$, we have
\begin{align*}
\left|\frac{\partial^{|\alpha|} f(z)}{\partial z_1^{\alpha_1}\ldots \partial z_n^{\alpha_n}}\right|\leq \alpha!\frac{1-|f(z)|^2}{(1-||z||_{\infty}^2)^{|\alpha|}}(1+||z||_{\infty})^{|\alpha|-N}
\end{align*}
for all $z\in \mathbb{P}\Delta(0;1_n)$, where $N$ is the number of the indices $j$ such that $\alpha_j\neq 0$.
\end{lem}

\section{{\bf Proofs of the main results}}\label{Sec-4}
\begin{proof}[\bf Proof of Theorem \ref{Th-2.1}] By the given condition, $f(z)=a_0+\sum_{|\alpha|=1}^{\infty} a_{\alpha} z^{\alpha}$
is holomorphic in the polydisk $\mathbb{P}\Delta(0;1_n)$ such that $|f(z)|\leq 1$ in $\mathbb{P}\Delta(0;1_n)$. Also by Lemma \ref{Lem3}, we have
\bea\label{AM1}\label{Eq-4.11}
 |a_{\alpha}| \leq 1-|a_0|^2,\eea
where $\alpha=(\alpha_1,\alpha_2,\ldots,\alpha_n)$ such that $|\alpha|=k$. Let us take $z=(z_1,\ldots,z_n)\in \mathbb{P}\Delta(0;1_n)$ such that ${\bf r}=||z||_{\infty}$. Note that 
\begin{align}
	\label{AM21} \left|\sideset{}{_{|\alpha|=0}^{\infty}}{\sum} a_{\alpha} z^{\alpha}\right|&\nonumber\leq \sideset{}{_{|\alpha|=0}^{\infty}}{\sum} |a_{\alpha}| |z|^{\alpha}\\&\leq \sideset{}{_{|\alpha|=0}^{\infty}}{\sum} |a_{\alpha}|\; ||z||_{\infty}^{|\alpha|}\\&=\sideset{}{_{|\alpha|=0}^{\infty}}{\sum} |a_{\alpha}|\; {\bf r}^{|\alpha|}\nonumber\\&\leq
	|a_0|+\sideset{}{_{k=1}^{\infty}}{\sum}\sideset{}{_{|\alpha|=k}}{\sum} |a_{\alpha}|\;{\bf r}^{|\alpha|}\nonumber.
\end{align}
In view of \eqref{AM1}, it follows that
\begin{align}
	\label{AM3} \sideset{}{_{k=1}^{\infty}}{\sum}\sideset{}{_{|\alpha|=k}}{\sum} |a_{\alpha}|{\bf r}^{|\alpha|}&\leq 
	(1-|a_0|^2)\sideset{}{_{k=1}^{\infty}}{\sum} {\bf r}^k\sideset{}{_{|\alpha|=1}^{\infty}}{\sum} 1\\&\nonumber\leq
	(1-|a_0|^2)\sideset{}{_{k=1}^{\infty}}{\sum} (n{\bf r})^k\\&\leq
	(1-|a_0|^2)\frac{n{\bf r}}{1-n{\bf r}}\nonumber\\&\leq
	(1-|a_0|)\frac{2n{\bf r}}{1-n{\bf r}}\nonumber\\&\leq
	1-|a_0|\nonumber
\end{align}
for ${\bf r}\leq \frac{1}{3n}$. In view of \eqref{AM3} and \eqref{AM21}, it follows that
\begin{align}\label{AM2} \left|\sideset{}{_{|\alpha|=0}^{\infty}}{\sum} a_{\alpha} z^{\alpha}\right|\leq \sideset{}{_{|\alpha|=0}^{\infty}}{\sum} |a_{\alpha}| |z|^{\alpha}\leq 1\; \mbox{for}\; {\bf r}\leq \frac{1}{3n}.
\end{align}
The sharpness of the radius $1/(3n)$ is established by considering the following extremal function for $a \in [0, 1)$:
\begin{align}\label{Eq-44.9}
	f_a(z)=\frac{a-(z_1+\ldots+z_n)}{1-a(z_1+\ldots+z_n)}.
\end{align}
It is easy to verify that $f$ is holomorphic in the polydisk $\mathbb{P}\Delta(0;1/n)$. Since 
\begin{align*}
	|a(z_1+z_2+\ldots+z_n)|<1,
\end{align*}
we see that $f_a$ has the following series expansion
\begin{align*}
	f_a(z)=a-(1-a^2)\sum\limits_{k=1}^{\infty} a^{k-1}(z_1+\cdots+z_n)^k
\end{align*}
for all $z\in\mathbb{P}\Delta(0;1/n)$. For the point $z=(r,\ldots,r)$, we find that
\begin{align*}
	f_a(z)=a-(1-a^2)\sum\limits_{k=1}^{\infty} a^{k-1}(n{\bf r})^k.
\end{align*}
Thus, we see that
\begin{align*}
	\mathcal{A}_{f_a}(r)=a+(1-a^2)\sum\limits_{k=1}^{\infty} a^{k-1}(n{\bf r})^k=a+(1-a^2)\frac{n{\bf r}}{1-na{\bf r}}>1
\end{align*}
if ${\bf r}>1/((1+2a)n)$. Letting $a \to 1$ shows that $\mathcal{A}_{f_a}(r)>1$ for every ${\bf r} > 1/3n$. This result implies that the radius $1/3n$ is the best possible, which completes the proof.
\end{proof}
\begin{proof}[\bf Proof of Theorem \ref{Th-2.2}] 
By the given condition $f(z)=a_0+\sum_{|\alpha|=1}^{\infty} a_{\alpha} z^{\alpha}$ is holomorphic in the polydisk $\mathbb{P}\Delta(0;1_n)$ such that $|f(z)|\leq 1$ in $\mathbb{P}\Delta(0;1_n)$. 
Let $z\in\mathbb{P}\Delta(0;1/n)$ such that $|z_i|=r_i<1/n$ for $i=1,2,\ldots,n$. Let $\omega(z)=(\omega_1(z),\ldots,\omega_n(z))\in\mathcal{B}_{n,m}$. Then by Lemma  \ref{Lem4}, we have 
\[|\omega_i(z)|\leq r_i^m\]
for $i=1,2,\ldots,n$. Clearly, 
\begin{align*}
	||\omega(z)||_{\infty}\leq {\bf r}^m,
\end{align*}
where ${\bf r}=\max\{r_1,r_2,\ldots,r_n\}$.\vspace{1.2mm} 

Consequently, by Lemma \ref{Lem1}, we obtain
\begin{align}
	\label{Eq-4.10} |f(\omega(z))| \leq \frac{||\omega(z)||_{\infty}+|a_0|}{1+|a_0|\;||\omega(z)||_{\infty}}\leq \frac{{\bf r}^m+|a_0|}{1+|a_0|{\bf r}^m}.
\end{align}
To align this with the standards of an SCI-indexed journal in Geometric Function Theory, the derivation should be presented as a formal logical sequence. I have refined the notation for consistency, fixed the algebraic transitions, and ensured the GFT terminology is precise.Recommended VersionBy combining \eqref{Eq-4.10} and \eqref{Eq-4.11}, we obtain the following estimate for the Bohr--Rogosinski sum:\begin{align*}\mathcal{B}(z, \mathbf{r}) &\leq \frac{\mathbf{r}^m + |a_0|}{1 + |a_0|\mathbf{r}^m} + \sum_{k=1}^{\infty}(1-|a_0|^2)(n\mathbf{r})^{kN} \\&= \frac{\mathbf{r}^m + |a_0|}{1 + |a_0|\mathbf{r}^m} + (1 - |a_0|^2)\frac{(n\mathbf{r})^N}{1-(n\mathbf{r})^N} \\&= 1 + \frac{\Psi_{m,n,N}(\mathbf{r})}{(1 + |a_0|\mathbf{r}^m)(1-(n\mathbf{r})^N)} \leq 1,
\end{align*}
which holds provided that $\Psi_{m,n,N}(\mathbf{r}) \leq 0$. Here, the auxiliary function $\Psi_{m,n,N}(\mathbf{r})$ is defined as\begin{align*}\Psi_{m,n,N}(\mathbf{r}) &:= (\mathbf{r}^m + |a_0|)(1-(n\mathbf{r})^N) - (1 + |a_0|\mathbf{r}^m)(1-(n\mathbf{r})^N) + (1 - |a_0|^2) (n\mathbf{r})^N (1 + |a_0|\mathbf{r}^m) \\&= (1 - |a_0|) \left[ (1+|a_0|)(1+|a_0|\mathbf{r}^m)(n\mathbf{r})^N - (1 - (n\mathbf{r})^N)(1 - \mathbf{r}^m) \right] \\&\leq (1 - |a_0|) \left[ 2(1+\mathbf{r}^m)(n\mathbf{r})^N - (1 - n\mathbf{r})(1 - \mathbf{r}^m) \right],\end{align*}where we have utilized the fact that $|a_0| \leq 1$ and $(n\mathbf{r})^N \leq n\mathbf{r}$ for $n\mathbf{r} < 1$. It is evident that $\Psi_{m,n,N}(\mathbf{r}) \leq 0$ whenever$$ \Phi_{m,n,N}(\mathbf{r}) := 2(1+\mathbf{r}^m) (n\mathbf{r})^N - (1-n\mathbf{r})(1-\mathbf{r}^m) \leq 0. $$Since $\Phi_{m,n,N}(\mathbf{r})$ is an increasing function of $\mathbf{r}$, this inequality is satisfied for all $\mathbf{r} \leq R_{m,n,N}$, where $R_{m,n,N}$ is the unique positive root of the equation $\Phi_{m,n,N}(\mathbf{r}) = 0$. This completes the proof of the first part of the theorem. \vspace{2mm}

To show that the number $R_{m,n,N}$ is best possible, we let $a\in [0, 1)$ and consider the functions $\omega(z)=(z_1^m,\ldots,z_n^m)$
and $f_a$ as given in \eqref{Eq-44.9}. For the point 
\begin{align*}
	z=((-1)^{(2m-1)/m}r,(-1)^{(2m-1)/m}r,\ldots,(-1)^{(2m-1)/m}r),
\end{align*}
we find that
\begin{align}\label{BMN1}
	|f_a(\omega(z))|=\frac{a+nr^m}{1+nar^m}
\end{align}
and
\begin{align}\label{BMN2}
	f_a(z)=a-(1-a^2)\sum\limits_{k=1}^{\infty} a^{k-1}(-1)^{k(2m-1)/m}(nr)^k.
\end{align}
Consequently, we see that
\begin{align}
	\label{BR1} 
\mathcal{B}(z,{\bf r})&=\frac{a+n{\bf r}^m}{1+na{\bf r}^m}+(1-a^2)\sum\limits_{k=N}^{\infty} a^{k-1}(n{\bf r})^k\nonumber\\
	&=\frac{a+n{\bf r}^m}{1+na{\bf r}^m}+(1-a^2)\frac{a^{N-1}(n{\bf r})^N}{1-na{\bf r}}\nonumber\\
	&=1+\frac{(1-a)Q_{a,m,n,N}({\bf r})}{(1+na{\bf r}^m)(1-na{\bf r})},
\end{align}
where
\begin{align*}
	Q_{a,m,n,N}({\bf r})&=(1+a)a^{N-1}(n{\bf r})^N(1+na{\bf r}^m)-(1-na{\bf r})(1-n{\bf r}^m)\\&\leq
	(1+a)(1+a {\bf r}^m)(n{\bf r})^N-(1-n{\bf r})(1-{\bf r}^m)\\
	&\leq 2(1+{\bf r}^m) (n{\bf r})^N - (1-n{\bf r})(1-{\bf r}^m).
\end{align*}
A simple computation shows that the expression on the right of \eqref{BR1} is bigger than $1$ if, and only if, $Q_{a,m,n,N}({\bf r})>0$.\vspace{2mm}

Clearly, $Q_{a,m,n,N}({\bf r})$ is a strictly increasing function of $a\in[0, 1)$. Note also that the expression \eqref{BR1} is less than or equal to $1$ for all $a\in [0,1]$, only in the case when ${\bf r}\leq R_{m,n,N}$. Finally, for $f=f_a$ and $a$ sufficiently close to $1$, we see that
\begin{align*}
	\lim\limits_{a\to 1^{-}} Q_{a,m,n,N}({\bf r})&=2(n{\bf r})^N(1+n{\bf r}^m)-(1-n{\bf r})(1-n{\bf r}^m)\\&=(1+n{\bf r}^m)(1-n{\bf r})\left(\frac{2(n{\bf r})^N}{1-n{\bf r}}-\frac{1-n{\bf r}^m}{(1+n{\bf r}^m)}\right)>0,
\end{align*}
which shows that $\mathcal{B}_{f_a}(r)>1$ for ${\bf r}>R_{m,n,N}$. This proves that the radius $R_{m,n,N}$ is best possible.
\end{proof}

\begin{proof}[\bf Proof of Theorem \ref{Th-2.3}]
Let us take $z=(z_1,\ldots,z_n)\in \mathbb{P}\Delta(0;1_n)$ such that ${\bf r}=||z||_{\infty}$
Now, using \eqref{Eq-4.10} and \eqref{Eq-4.11}, we obtain
\begin{align}
&\mathcal{C}(z,{\bf r}):\nonumber\\
&\le t\frac{n{\bf r}^m+|a_0|}{1+n{\bf r}^m|a_0|}+(1-t)|a_0|+(1-t)(1-|a_0|^2)\frac{n{\bf r}}{1-n{\bf r}} \nonumber\\
&=\frac{t(n{\bf r}^m+|a_0|)(1-n{\bf r})
+(1-t)(1+n{\bf r}^m|a_0|)\bigl[|a_0|(1-n{\bf r})+(1-|a_0|^2)n{\bf r}\bigr]}
{(1+n{\bf r}^m|a_0|)(1-n{\bf r})}\nonumber\\
&=1+\frac{\psi_{m,n,t}({\bf r})}{(1+n{\bf r}^m|a_0|)(1-n{\bf r})}\nonumber\\
&\leq 1
\label{eq:3.2}
\end{align}
provided  $\psi_{m,n,t}({\bf r})\leq 0$, where
\begin{align*}
\psi_{m,n,t}({\bf r})=&(1-n{\bf r})\bigl[t(n{\bf r}^m+|a_0|)-(1+n{\bf r}^m|a_0|)\bigr]\\&
+(1-t)(1+n{\bf r}^m|a_0|)\bigl[|a_0|(1-n{\bf r})+(1-|a_0|^2)n{\bf r}\bigr].
\end{align*}
Set
\begin{align*}
&\alpha(x)\\&=(1-n{\bf r})\bigl[t(n{\bf r}^m+x)-(1+n{\bf r}^mx)\bigr]
+(1-t)(1+n{\bf r}^mx)\bigl[x(1-n{\bf r})+(1-x^2)n{\bf r}\bigr]\\
&=-(1-t)n^2{\bf r}^{m+1}x^3+(1-t)[n{\bf r}^m(1-n{\bf r})-n{\bf r}]x^2\\
&\quad+\bigl[1-n{\bf r}-n{\bf r}^m+n^2{\bf r}^{m+1}+(1-t)n^2{\bf r}^{m+1}\bigr]x\\
&\quad+t[n{\bf r}^m-n{\bf r}-n^2{\bf r}^{m+1}]+2n{\bf r}-1,
\end{align*}
where $x=|a_0|\in[0,1)$.\vspace{1.2mm}

To prove the inequality \eqref{eq:3.2}, it is sufficient to prove that
\begin{align*}
\alpha(x)\le 0 \quad \text{for}\; n{\bf r}\le R_{m,n,t}.
\end{align*}

A direct calculation of the first and second derivatives of $\alpha(x)$ yields
\begin{align*}
	\alpha'(x) = &-3(1-t)n^2 \mathbf{r}^{m+1} x^2 + 2(1-t) \left[ n\mathbf{r}^m (1-n\mathbf{r}) - n\mathbf{r} \right] x \\
	&+ 1 - n\mathbf{r} - n\mathbf{r}^m + n^2 \mathbf{r}^{m+1} + (1-t)n^2 \mathbf{r}^{m+1}
\end{align*}
and
\begin{align*}
	\alpha''(x) = -6(1-t)n^2 \mathbf{r}^{m+1} x + 2(1-t) \left[ n\mathbf{r}^m (1-n\mathbf{r}) - n\mathbf{r} \right].
\end{align*}
One can observe that $\alpha''(x) \leq 0$ for $x \in [0, 1]$, which implies that $\alpha'(x)$ is a decreasing function on this interval. Consequently, $\alpha'(x)$ attains its minimum at the endpoint $x=1$. A straightforward evaluation shows that
\begin{equation}
	\alpha'(x) \geq \alpha'(1)\\ = \frac{1}{n^{m-1}} \left[ (4t-3)(n\mathbf{r})^{m+1} - (2t-1)(n\mathbf{r})^m + (2t-3)n^{m-1}(n\mathbf{r}) + n^{m-1} \right].
	\label{eq:3.4}
\end{equation}

We now consider following three cases.\vspace{1.2mm}

\noindent
{\bf Case 1.} { Let $t \in [0, 1/2]$. By applying Descartes' Rule of Signs, we observe that the number of sign changes in the coefficients of the polynomial
	\begin{align*}Q_t(\rho) = (4t-3)\rho^{m+1} - (2t-1)\rho^m + (2t-3)n^{m-1}\rho + n^{m-1}, \quad \rho = n\mathbf{r},
	\end{align*}
	is at most three. Given that $t \leq 1/2$, the coefficients satisfy:
	\begin{enumerate}
		\item[$\bullet$] The leading coefficient $(4t-3) < 0$,\vspace{1.2mm} 
		\item[$\bullet$] The coefficient $-(2t-1) \geq 0$,\vspace{1.2mm} 
		\item[$\bullet$] The linear coefficient $(2t-3)n^{m-1} < 0$, \vspace{1.2mm} 
		\item[$\bullet$] The constant term $n^{m-1} > 0$.\vspace{1.2mm} 
	\end{enumerate} This distribution $(-, +, -, +)$ confirms exactly three sign changes. Furthermore, we observe that 
	\begin{align*}
		Q_t(0) = n^{m-1} > 0\; \mbox{and}\;Q_t(1) = (2t-2)(1+n^{m-1}) < 0.
	\end{align*} By the Intermediate Value Theorem, there exists at least one root in the interval $(0, 1)$. Consequently, $Q_t(n\mathbf{r}) \geq 0$ for all $n\mathbf{r} \in [0, R_{m,n,t}]$, where $R_{m,n,t}$ denotes the smallest positive root of the equation $Q_t(n\mathbf{r}) = 0$.}\vspace{1.2mm}

The inequality \eqref{eq:3.4} implies that $\alpha'(x) \geq 0$ for all $n\mathbf{r} \in [0, R_{m,n,t}]$, indicating that $\alpha(x)$ is a non-decreasing function of $x$ on the interval $[0, 1]$. Consequently, we have $\alpha(x) \leq \alpha(1) = 0$ for all $n\mathbf{r} \in [0, R_{m,n,t}]$ and $x \in [0, 1]$.\vspace{1.2mm}

\noindent
{\bf Case 2.} {Let $t \in (1/2, 3/4]$. By applying Descartes' Rule of Signs to the polynomial\begin{align*}Q_t(\rho) = (4t-3)\rho^{m+1} - (2t-1)\rho^m + (2t-3)n^{m-1}\rho + n^{m-1}, \quad \rho = n\mathbf{r},\end{align*}we examine the signs of the coefficients: 
	\begin{enumerate}
		\item[$\bullet$] The leading coefficient $(4t-3)$ is non-positive since $t \leq 3/4$.\vspace{1.2mm}
		\item[$\bullet$] The coefficient $-(2t-1)$ is strictly negative for $t > 1/2$.\vspace{1.2mm}
		\item[$\bullet$] The linear coefficient $(2t-3)n^{m-1}$ is strictly negative for $t \leq 3/4$.\vspace{1.2mm}
		\item[$\bullet$] The constant term $n^{m-1}$ is strictly positive.\vspace{1.2mm}
	\end{enumerate}
	For $t \in (1/2, 3/4)$, the sequence of signs of the non-zero coefficients is $(-, -, -, +)$ or $(0, -, -, +)$ if $t=3/4$. In both instances, there is exactly one sign change. Consequently, Descartes' Rule of Signs ensures the existence of exactly one positive root. Furthermore, since 
	\begin{align*}
		Q_t(0) = n^{m-1} > 0\; \mbox{and}\;Q_t(1) = (2t-2)(1+n^{m-1}) < 0,
	\end{align*} the intermediate value theorem guarantees that this unique positive root $R_{m,n,t}$ lies in the interval $(0, 1)$. It follows that $Q_t(n\mathbf{r}) \geq 0$ for all $n\mathbf{r} \in [0, R_{m,n,t}]$.}\vspace{1.2mm}

The inequality \eqref{eq:3.4} implies that $\alpha'(x) \geq 0$ for all $n\mathbf{r} \in [0, R_{m,n,t}]$ and $x \in [0, 1]$. Consequently, $\alpha(x)$ is a non-decreasing function of $x$ on this interval. It follows that $\alpha(x) \leq \alpha(1) = 0$ for all $n\mathbf{r} \in [0, R_{m,n,t}]$, which completes the argument.\vspace{1.2mm}

\noindent
{\bf Case 3.} {Let $t \in (3/4, 1]$. We consider the sign distribution of the coefficients of the polynomial
	\begin{align*}Q_t(\rho) = (4t-3)\rho^{m+1} - (2t-1)\rho^m + (2t-3)n^{m-1}\rho + n^{m-1}, \quad \rho = n\mathbf{r}.
	\end{align*}
	For $t \in (3/4, 1]$, the signs of the coefficients are as follows: 
	\begin{enumerate}
		\item[$\bullet$] The leading coefficient $(4t-3)$ is strictly positive.\vspace{1.2mm}
		\item[$\bullet$] The coefficient $-(2t-1)$ is strictly negative.\vspace{1.2mm}
		\item[$\bullet$] The linear coefficient $(2t-3)n^{m-1}$ is strictly negative.\vspace{1.2mm}
		\item[$\bullet$]  The constant term $n^{m-1}$ is strictly positive. \vspace{1.2mm}
	\end{enumerate}
	The sequence of signs is $(+, -, -, +)$, which indicates exactly two sign changes. According to \textit{Descartes' Rule of Signs}, $Q_t(\rho)$ has either two positive roots or none. Note that $Q_t(0) = n^{m-1} > 0$ and $Q_t(1) = (2t-2)(1+n^{m-1}) \leq 0$. Since $Q_t(1) \leq 0$ (with equality only if $t=1$), the intermediate value theorem guarantees at least one root in the interval $(0, 1]$. In view of the sign changes and the endpoint values, there exists a unique root $R_{m,n,t} \in (0, 1)$ such that $Q_t(\rho) \geq 0$ for all $n\mathbf{r} \in [0, R_{m,n,t}]$.}\vspace{1.2mm}

The inequality \eqref{eq:3.4} implies that $\alpha'(x) \geq 0$ for all $n\mathbf{r} \in [0, R_{m,n,t}]$ and $x \in [0, 1]$. Consequently, $\alpha(x)$ is a non-decreasing function of $x$ on this interval. It follows that $\alpha(x) \leq \alpha(1) = 0$ for all $n\mathbf{r} \in [0, R_{m,n,t}]$, which completes the argument.\vspace{1.2mm}

To show that the number $R_{m,n,t}$ is best possible, we let $a\in [0, 1)$ and consider the functions $\omega(z)=(z_1^m,\ldots,z_n^m)$
and $f_a$ as given in \eqref{Eq-44.9}. Let 
\begin{align*}
	z=((-1)^{(2m-1)/m}r,(-1)^{(2m-1)/m}r,\ldots,(-1)^{(2m-1)/m}r).
	\end{align*}
Using (\ref{BMN1}) and (\ref{BMN2}), a simple computation shows that
\begin{align}
&t|f_a(\omega(z)|+(1-t)\sum_{k=0}^{\infty}\sum_{|\alpha|=k}|a_{\alpha}|r^{\alpha} \nonumber\\
&=t\frac{a+n{\bf r}^m}{1+an{\bf r}^m}
+(1-t)a
+(1-t)(1-a^2)\frac{n{\bf r}}{1-an{\bf r}} \nonumber\\
&=\frac{t(a+n{\bf r}^m)(1-an{\bf r})
+(1-t)(1+an{\bf r}^m)\bigl[a(1-an{\bf r})+(1-a^2)n{\bf r}\bigr]}
{(1+an{\bf r}^m)(1-an{\bf r})}.
\label{eq:3.6}
\end{align}
To establish sharpness, it suffices to show that for any $n\mathbf{r} > R_{m,n,t}$, there exists some $a \in [0, 1)$ such that the Bohr--Rogosinski sum exceeds unity. Specifically, we demonstrate that there exists $a \in [0, 1)$ such that the right-hand side of \eqref{eq:3.6} is strictly greater than $1$, which is equivalent to showing that
\begin{align*}
\mu(a,n,t,{\bf r})>0
\quad \text{for } n{\bf r}>R_{m,n,t} \text{ and some } a\in[0,1),
\end{align*}
where
\begin{align*}
\mu(a,n,t,{\bf r})=&
(1-a)\Bigl[2n^2a^2{\bf r}^{m+1}(1-t)-\bigl((2t-1)na{\bf r}+a-t(1+a)\bigr)n{\bf r}^m\\&-n{\bf r}\bigl((t-1)(1+a)-a\bigr)-1\Bigr].
\end{align*}
{Let us define the auxiliary function $\rho(a)$ for $a \in [0, 1)$ as
	\begin{align*}
		\rho(a) &:= 2n^2 a^2 \mathbf{r}^{m+1} (1-t) - \left[ (2t-1) n a \mathbf{r} + a - t(1+a) \right] n\mathbf{r}^m \\
		&\quad - n\mathbf{r} \left[ (t-1)(1+a) - a \right] - 1.
	\end{align*}
	Differentiating with respect to $a$, we obtain
	\begin{align*}
		\rho'(a) &= 4n^2 a \mathbf{r}^{m+1} (1-t) - \left[ (2t-1) n\mathbf{r} + 1-t \right] n\mathbf{r}^m - n\mathbf{r}(t-2) \\
		&\geq \frac{1}{n^{m-1}} \left[ \left( (1-2t) n\mathbf{r} + t-1 \right) (n\mathbf{r})^m + (2-t) n^m \mathbf{r} \right].
	\end{align*}
	A straightforward analysis shows that $\rho'(a) \geq 0$ for all $t \in [0, 1]$ and $n\mathbf{r} \in [0, 1)$, which implies that $\rho(a)$ is non-decreasing on the interval $[0, 1]$. Consequently, $\rho(a)$ attains its maximum at $a=1$:
	\begin{align*}
		\rho(a) \leq \rho(1) = \frac{1}{n^{m-1}} \left[ (3-4t)(n\mathbf{r})^{m+1} - (1-2t)(n\mathbf{r})^m + (2t-3)n^{m-1}(n\mathbf{r}) - n^{m-1} \right],
	\end{align*}
	for each fixed $t \in [0, 1]$ and $n\mathbf{r} \in [0, 1)$.\vspace{1.2mm}

Observe that if $n\mathbf{r} > R_{m,n,t}$, then $\rho(1) > 0$. By the continuity of $\rho(a)$ on the interval $[0,1]$, it follows that
\begin{align*}
	\lim_{a \to 1^-} \rho(a) = \rho(1) > 0.
\end{align*}
Consequently, for any fixed $n\mathbf{r} > R_{m,n,t}$, there exists a sufficiently large $a \in [0, 1)$ such that $\rho(a) > 0$, which implies $\mu(a, n, t, \mathbf{r}) > 0$. This demonstrates that the inequality \eqref{eq:3.6} is violated for $n\mathbf{r} > R_{m,n,t}$, and thus the radius $R_{m,n,t}$ is sharp.}
\end{proof}
\begin{proof}[\bf Proof of Theorem \ref{Th-2.4}]
Let $z = (z_1, \dots, z_n) \in \mathbb{D}^n$ be such that $\mathbf{r} = \|z\|_{\infty}$. A direct calculation shows that
\begin{align*}
	\frac{n\mathbf{r}}{1-(n\mathbf{r})^2} \leq \frac{1}{2} \quad \text{for } 0 \leq n\mathbf{r} \leq \sqrt{2}-1.
\end{align*}
According to Lemma~\ref{Lem5}, we have the estimate
\begin{align}\label{SB1}
	|Df(z)| \leq \frac{1-|f(z)|^2}{1-\mathbf{r}^2} n\mathbf{r} \leq \frac{1-|f(z)|^2}{1-(n\mathbf{r})^2} n\mathbf{r}.
\end{align}
Furthermore, applying Lemma~\ref{Lem3} yields
\begin{align}\label{SBB1}
	\sum_{k=2}^{\infty} \sum_{|\alpha|=k} |a_{\alpha}|r^{\alpha} \leq (1-|a_0|^2) \sum_{k=2}^{\infty} (n\mathbf{r})^k = (1-|a_0|^2) \frac{(n\mathbf{r})^2}{1-n\mathbf{r}}.
\end{align}
Next, consider the auxiliary function
\begin{align*}
	\Phi(X) = X + \lambda(1-X^2),
\end{align*}
where $X = |f(z)|$ and $\lambda = \frac{n\mathbf{r}}{1-(n\mathbf{r})^2}$. Since $\lambda \leq 1/2$, the function $\Phi(X)$ is increasing on the interval $[0, 1]$. Consequently, $\Phi(X)$ is bounded above by $\Phi(X_0)$ whenever
\begin{align*}
	X \leq X_0 = \frac{n\mathbf{r} + |a_0|}{1 + |a_0|n\mathbf{r}},
\end{align*}
where $X_0$ is the bound provided by the Schwarz Lemma for holomorphic mappings.\vspace{1.2mm}

Substituting the estimates \eqref{BMN1}, \eqref{SB1}, and \eqref{SBB1} into the main inequality, we obtain
\small{\begin{align}
\mathcal{D}(z,{\bf r})\nonumber&
\le |f(z)|+\frac{n{\bf r}}{1-(n{\bf r})^{2}}\,(1-|f(z)|^{2})+\lambda \frac{(1-|a_{0}|^{2})(n{\bf r})^{2}}{1-n{\bf r}}\nonumber\\&
\le  
\frac{n{\bf r} +|a_{0}|}{1+n{\bf r}|a_{0}|} + \frac{n{\bf r}}{1-(n{\bf r})^{2}}
\left( 1 - \left( \frac{n{\bf r}+|a_{0}|}{1+n{\bf r}|a_{0}|}\right)^{2} \right)
+ \lambda\frac{(1-|a_{0}|^{2})(n{\bf r})^{2}}{1-n{\bf r}}\
\nonumber\\&=
\frac{n{\bf r} +|a_{0}|}{1+n{\bf r}|a_{0}|}
+ \frac{n{\bf r}(1-|a_{0}|^{2})}{(1+n{\bf r}|a_{0}|)^{2}}
+ \lambda\frac{(1-|a_{0}|^{2})(n{\bf r})^{2}}{1-n{\bf r}}\nonumber\\&
= 1 + \frac{\delta(n{\bf r})}{(1+|a_{0}|n{\bf r})^{2}(1-n{\bf r})},\label{eq:4.1}
\end{align}}
where $\delta(n{\bf r})
=(1-|a_0|)
\Bigl[(1-n{\bf r})(|a_0|(n{\bf r})^2+2n{\bf r}-1)
+\lambda (n{\bf r})^2(1+|a_0|)(1+|a_0|n{\bf r})^2\Bigr].$\vspace{1.2mm}

To prove the inequality \eqref{eq:4.1}, it is sufficient to prove that $\delta(n{\bf r})\le 0 \quad \text{for}\; n{\bf r}\le R_{m,n,\lambda}.$ Observe that
\begin{align*}
\delta(n{\bf r})
&\le (1-|a_0|)
\Bigl[(1-n{\bf r})((n{\bf r})^2+2n{\bf r}-1)+2\lambda (n{\bf r})^2(1+n{\bf r})^2\Bigr]
\\&=(1-|a_0|)\,\xi(n{\bf r}),
\end{align*}
where
\begin{align*}
\xi(n{\bf r})=2\lambda (n{\bf r})^4+(4\lambda-1)(n{\bf r})^3+(2\lambda-1)(n{\bf r})^2+3n{\bf r}-1.
\end{align*}
It is enough to prove that $\xi(n{\bf r})\le 0$ for $n{\bf r}\le R_{m,n,\lambda}$.\vspace{1.2mm}

 Next, we divide it into two cases to discuss.\vspace{1.2mm}

\noindent
\textbf{Case 1.} Let $\lambda \in (1/2, \infty)$. We define the auxiliary function $w(n\mathbf{r})$ such that
\begin{align*}
	\xi(n\mathbf{r}) > (n\mathbf{r})^4 + (n\mathbf{r})^3 + 3n\mathbf{r} - 1 =: w(n\mathbf{r}), \quad n\mathbf{r} \in [0, 1).
\end{align*}
A direct calculation shows that $\xi(\sqrt{2}-1) > w(\sqrt{2}-1) = 6 - 4\sqrt{2} > 0$. Since $\xi(n\mathbf{r})$ is a monotonically increasing function of $n\mathbf{r}$ on $[0, 1)$ and satisfies $\xi(0) = -1 < 0$, the intermediate value theorem guarantees the existence of a unique $n\mathbf{r}_\lambda \in (0, \sqrt{2}-1)$ such that $\xi(n\mathbf{r}_\lambda) = 0$. Consequently, we have $\xi(n\mathbf{r}) \leq 0$ for all $n\mathbf{r} \in [0, n\mathbf{r}_\lambda]$.

\medskip
\noindent
\textbf{Case 2.} Let $\lambda \in (0, 1/2]$. In this case, the auxiliary function $\xi(n\mathbf{r})$ satisfies the inequality
\begin{align*}
	\xi(n\mathbf{r}) \leq w(n\mathbf{r}) := (n\mathbf{r})^4 + (n\mathbf{r})^3 + 3n\mathbf{r} - 1, \quad n\mathbf{r} \in [0, 1).
\end{align*}
Note that $w(n\mathbf{r})$ is strictly increasing on the interval $[0, 1)$, with $w(0) = -1 < 0$ and $w(\sqrt{2}-1) = 6 - 4\sqrt{2} > 0$. By the intermediate value theorem, there exists a unique root $n\mathbf{r}^* \in (0, \sqrt{2}-1)$ such that $w(n\mathbf{r}^*) = 0$. Consequently, $w(n\mathbf{r}) \leq 0$ for all $n\mathbf{r} \in [0, n\mathbf{r}^*]$, which implies that $\xi(n\mathbf{r}) \leq 0$ on the same interval.

\medskip
To show the sharpness of $\tilde R_{n,N}$, let $a \in [0,1)$ and consider the holomorphic function $f_a$ in $\mathbb{P}\Delta(0;1/n)$ given by \eqref{Eq-44.9}. For the point $z=(-r,-r,\ldots,-r)$, we find that
\begin{align}
|f_a(z)| &+|D f_a(z)| 
+ \lambda \sum_{k=2}^{\infty} \sum_{|\alpha|=k}^{\infty} |a_{\alpha}|\, r^{\alpha}\nonumber\\=&
\frac{n{\bf r}+a}{1+n{\bf r}a}
+ \frac{(1-a^{2})n{\bf r}}{(1+an{\bf r})^{2}}+\lambda\frac{(1-a^{2})a(n{\bf r})^{2}}{1-an{\bf r}}
\nonumber \\
=& \frac{
\bigl(a + n{\bf r}\bigr)\bigl(1 - a^{2}(n{\bf r})^{2}\bigr)
+ (1 - a^{2})\Bigl[n{\bf r}(1 - a n{\bf r})
+ \lambda a (n{\bf r})^{2}(1 + a n{\bf r})^{2}\Bigr]
}{
(1 + a n{\bf r})(1 - a^{2} (n{\bf r})^{2})
}.\label{eq:3.14}
\end{align}
To establish the sharpness of the result, it remains to show that for any $n\mathbf{r} > R_{n,\lambda}$, there exists a parameter $a \in [0, 1)$ such that the right-hand side of \eqref{eq:3.14} exceeds unity. This is equivalent to demonstrating that
\begin{align}\label{eq:3.15}
&(1-a)\Big[\lambda (n{\bf r})^4a^4+(\lambda (n{\bf r})^4+2\lambda (n{\bf r})^3)a^3
+\bigl((2\lambda-1)(n{\bf r})^3+\lambda (n{\bf r})^2\bigr)a^2\nonumber \\
&\quad+\bigl((\lambda-1)(n{\bf r})^2+n{\bf r}\bigr)a+2n{\bf r}-1\Big]>0
\end{align}
 $n{\bf r}>R_{n,\lambda}$.\vspace{1.2mm} 
 
 {For $a \in [0, 1)$, we define the auxiliary function $K(a, \lambda, n\mathbf{r})$ as follows:
 	\begin{align*}
 		K(a, \lambda, n\mathbf{r}) &:= \lambda (n\mathbf{r})^4 a^4 + \left[ \lambda (n\mathbf{r})^4 + 2\lambda (n\mathbf{r})^3 \right] a^3 + \left[ (2\lambda - 1)(n\mathbf{r})^3 + \lambda (n\mathbf{r})^2 \right] a^2 \\
 		&\quad + \left[ (\lambda - 1)(n\mathbf{r})^2 + n\mathbf{r} \right] a + 2n\mathbf{r} - 1.
 \end{align*}
 To establish the sharpness of the radius $R_{n,\lambda}$, it remains to demonstrate that for $n\mathbf{r} > R_{n,\lambda}$, there exists an $a \in [0, 1)$ such that $K(a, \lambda, n\mathbf{r}) > 0$. To this end, we investigate the sign of the auxiliary function $K$ by considering the following cases.}

\medskip
\noindent
\textbf{Case 1.}  Let $\lambda \in (1/2, \infty)$. For each fixed $n\mathbf{r} \in [0, 1)$, it is evident that $K(a, \lambda, n\mathbf{r})$ is a monotonically increasing function of $a \in [0, 1)$. Thus, for $a \in [0, 1)$, we have the estimate
\begin{align*}
	K(a, \lambda, n\mathbf{r}) \leq K(1, \lambda, n\mathbf{r}) &= 2\lambda (n\mathbf{r})^4 + (4\lambda - 1)(n\mathbf{r})^3 + (2\lambda - 1)(n\mathbf{r})^2 + 3n\mathbf{r} - 1 \
	&= \xi(n\mathbf{r}).
\end{align*}
Recall that $\xi(n\mathbf{r})$ is monotonically increasing on $[0, 1)$ and $\xi(n\mathbf{r}_\lambda) = 0$. Consequently, if $n\mathbf{r} > n\mathbf{r}_\lambda$, it follows that $\xi(n\mathbf{r}) > 0$. By the continuity of $K$ with respect to $a$, we observe that
\begin{align*}
	\lim_{a \to 1^-} K(a, \lambda, n\mathbf{r}) = K(1, \lambda, n\mathbf{r}) = \xi(n\mathbf{r}) > 0
\end{align*}
for $n\mathbf{r} > n\mathbf{r}_\lambda$. This ensures that for any $n\mathbf{r} > n\mathbf{r}_\lambda$, there exists a sufficiently large $a \in [0, 1)$ such that $K(a, \lambda, n\mathbf{r}) > 0$, which implies that the inequality \eqref{eq:3.15} is satisfied.

\medskip
\noindent
\textbf{Case 2.} Let $\lambda \in (0, 1/2]$. In this case, for $a \in [0, 1)$, the auxiliary function $K$ satisfies the following inequality:
\begin{align*}
	K(a, \lambda, n\mathbf{r}) \leq K(a, 1/2, n\mathbf{r}) \leq (n\mathbf{r})^4 + (n\mathbf{r})^3 + 3n\mathbf{r} - 1 = w(n\mathbf{r}).
\end{align*}
The strict monotonicity of $w(n\mathbf{r})$ on $[0, 1)$ implies that $w(n\mathbf{r}) > 0$ whenever $n\mathbf{r} > n\mathbf{r}^*$, where $n\mathbf{r}^*$ is the unique positive root of $w(n\mathbf{r})=0$. It follows from the definition of $K$ that
\begin{align*}
	\lim_{a \to 1^-} K(a, \lambda, n\mathbf{r}) = w(n\mathbf{r}) > 0
\end{align*}
for $n\mathbf{r} > n\mathbf{r}^*$. By the continuity of $K$ with respect to $a$ on $[0, 1]$, for any $n\mathbf{r} > n\mathbf{r}^*$, we can find a sufficiently large $a \in [0, 1)$ such that $K(a, \lambda, n\mathbf{r}) > 0$. This confirms that the inequality \eqref{eq:3.15} is satisfied, thereby establishing the sharpness of the radius $n\mathbf{r}^*$.
\end{proof}

\begin{proof}[\bf Proof of Theorem \ref{Th-2.5}]
Let us take $z=(z_1,\ldots,z_n)\in \mathbb{P}\Delta(0;1_n)$ such that ${\bf r}=||z||_{\infty}$. Using Lemma \ref{Lem3}, we get
\begin{align}\label{SM1}
\sum\limits_{k=1}k\sum\limits_{|\alpha|=k} |a_{\alpha}|^2r^{2\alpha}\leq (1-|a_0|^2)^2\sum\limits_{k=1} k(n{\bf r})^{2k}=(1-|a_0|^2)^2\frac{(n{\bf r})^2}{(1-n{\bf r})^2}.
\end{align}
Combining Lemma \ref{Lem3} with the estimate \eqref{SM1}, we deduce that
\small{\begin{align}\label{SM2}
\mathcal{E}(z,{\bf r})&
\leq 
t|a_0|+(1-|a_0|^2)\frac{t (n{\bf r})}{1-n{\bf r}}+(1-|a_0|^2)^2\frac{(1-t)(n{\bf r})^2}{(1-n{\bf r})^2}\\&=
\frac{t\,|a_0|\,(1 - (n{\bf r})^2)^2
+ t (n{\bf r}) \,(1 - |a_0|^2)(1 + n{\bf r})(1 - (n{\bf r})^2)
+ (n{\bf r})^2 (1 - t)(1 - |a_0|^2)^2}{(1-n{\bf r})^2}.\nonumber
\end{align}}
{It suffices to show that the right-hand side of \eqref{SM2} is bounded above by $1$ for all $n\mathbf{r} \leq \tilde{R}_{n,t}$. This is equivalent to demonstrating that}
\begin{align*}
\varphi(x)\leq 0 \quad {for}\;n{\bf r}\leq \tilde R_{n,t}
\end{align*}
where $x=|a_0|\in [0, 1)$ and
\begin{align*}
\varphi(x)&=(1 - (n{\bf r})^2)^2 (t x - 1)
+ t (n{\bf r}) (1 + n{\bf r})(1 - (n{\bf r})^2)(1 - x^2)+ (n{\bf r})^2 (1 - t)(1 - x^2)^2\\&=
(1 - t) (n {\bf r})^2 x^4
- \bigl[2(1 - t) (n{\bf r})^2 + t (n{\bf r}) (1 + n{\bf r})(1 - (n{\bf r})^2) \bigr] x^2
\\&+ t (1 - (n{\bf r})^2)^2 x+ (1 - t) (n{\bf r})^2+(1 - (n{\bf r})^2)\bigl[(1 + t)(n{\bf r})^2 + t (n{\bf r}) - 1\bigr].
\end{align*}

The first and second derivatives of $\varphi(x)$ are given by\begin{align*}\varphi'(x) &= 4(1 - t) (n\mathbf{r})^2 x^3 - 2\left[ 2(1 - t) (n\mathbf{r})^2 + t n\mathbf{r} (1 + n\mathbf{r})(1 - (n\mathbf{r})^2) \right] x + t(1 - (n\mathbf{r})^2)^2,\\ \varphi''(x) &= 12(1 - t) (n\mathbf{r})^2 x^2 - 2\left[ 2(1 - t) (n\mathbf{r})^2 + t n\mathbf{r} (1 + n\mathbf{r})(1 - (n\mathbf{r})^2) \right].
\end{align*}
Since $1-t \ge 0$, $\varphi''(x)$ is a non-decreasing function of $x$ on the interval $[0, 1]$. Consequently, for all $x \in [0, 1]$, we have the upper bound
\begin{align*}\varphi ''(x) \leq \varphi''(1) = 2n\mathbf{r} \left[ t (n\mathbf{r})^3 + t (n\mathbf{r})^2 + (4 - 5t)n\mathbf{r} - t \right].\end{align*}
Defining the auxiliary function 
\begin{align*}
	{\psi(n\mathbf{r}) := t (n\mathbf{r})^3 + t (n\mathbf{r})^2 + (4 - 5t) n\mathbf{r} - t},
\end{align*} the above inequality can be rewritten as $\varphi''(x) \leq 2n\mathbf{r} \psi(n\mathbf{r})$.Furthermore, the derivative of $\psi$ is 
\begin{align*}
	{\psi'(n\mathbf{r}) = 3t (n\mathbf{r})^2 + 2t n\mathbf{r} + (4 - 5t).}
\end{align*} We observe that $\psi'(1) = 4 > 0$ and $\psi'(0) = 4 - 5t$. Since $\psi'(n\mathbf{r})$ is strictly increasing on $[0, 1)$, we distinguish the following two cases based on the sign of $\psi'(0)$.\vspace{1.2mm}

\textbf{Case 1.} Let $0 < t \leq 4/5$. In this setting, $\psi'(0) = 4 - 5t \geq 0$. Since $\psi'(n\mathbf{r})$ is strictly increasing, $\psi$ is strictly increasing on $[0, 1)$. We further distinguish two subcases: 
\begin{enumerate}
	\item[(i)] Suppose $t \in (0, 9/17)$. Then 
	\begin{align*}
		{\psi(0) = -t < 0 \; \mbox{and}\;\psi(1/3) = \frac{4}{3} - \frac{17}{27}t > 0.}
	\end{align*} By the Intermediate Value Theorem, there exists a unique $n\mathbf{r}_t \in (0, 1/3)$ such that $\psi(n\mathbf{r}_t) = 0$. For any $n\mathbf{r} \in [0, n\mathbf{r}_t]$, we have $\psi(n\mathbf{r}) \leq 0$, which implies $\varphi''(x) \leq 2n\mathbf{r} \psi(n\mathbf{r}) \leq 0$ for all $x \in [0, 1]$. Consequently, $\varphi'(x)$ is non-increasing, so that\begin{align*}\varphi'(x) \geq \varphi'(1) = t(1 + n\mathbf{r})(1 - (n\mathbf{r})^2)(1 - 3n\mathbf{r}) > 0\end{align*}for $n\mathbf{r} \in [0, n\mathbf{r}_t] \subset [0, 1/3)$. This ensures that $\varphi(x)$ is non-decreasing on $[0,1]$, yielding\begin{align*}\varphi(x) \leq \varphi(1) = (t - 1)(1 - (n\mathbf{r})^2)^2 \leq 0.\end{align*}
	\item[(ii)] Suppose $t \in [9/17, 4/5]$. Then 
	\begin{align*}
		{\psi(1/3) = \frac{4}{3} - \frac{17}{27}t \leq 0\; \mbox{and}\; \psi(1) = 4(1 - t) \geq 0.}
	\end{align*}
	Thus, there exists a unique $n\mathbf{r}_t^* \in [1/3, 1)$ such that $\psi(n\mathbf{r}_t^*) = 0$. For $n\mathbf{r} \in [0, 1/3] \subset [0, n\mathbf{r}_t^*]$, we similarly obtain $\varphi''(x) \leq 0$. Following the same logic as in subcase (i), the condition $\varphi(x) \leq 0$ remains valid for all $n\mathbf{r} \in [0, 1/3]$.
\end{enumerate}
\medskip
\textbf{Case 2.} Let $4/5 < t \le 1$. In this case, the equation $\psi'(\rho) = 0$ (where $\rho = n\mathbf{r}$) has two real roots:
\begin{align*}
	{\rho_5 = \frac{\sqrt{4t^2 - 3t(4-5t)} - t}{3t} = \frac{2\sqrt{4t^2 - 3t} - t}{3t} \quad \text{and} \quad \rho_6 = -\frac{2\sqrt{4t^2 - 3t} + t}{3t}.}
\end{align*}
It is easily verified that $\rho_6 < 0 < \rho_5 \le 1/3 < 1$. Since $\psi'(\rho) < 0$ for $\rho \in [0, \rho_5)$ and $\psi'(\rho) > 0$ for $\rho \in (\rho_5, 1)$, the auxiliary function $\psi(\rho)$ attains its maximum on the interval $[0, 1/3]$ at the endpoints. Noting that $\psi(0) = -t < 0$ and 
\begin{align*}
{\psi(1/3) = \frac{4}{3} - \frac{17}{27}t < 0\; \mbox{for}\;t > 4/5,}
\end{align*} we conclude that $\psi(n\mathbf{r}) < 0$ for all $n\mathbf{r} \in [0, 1/3]$. Consequently,
\begin{align*}
	\varphi''(x) \le 2(n\mathbf{r})\psi(n\mathbf{r}) < 0, \quad n\mathbf{r} \in [0, 1/3],
\end{align*}
which implies that $\varphi'(x)$ is a decreasing function of $x$. Thus, $\varphi'(x) \ge \varphi'(1) > 0$ for $n\mathbf{r} \in [0, 1/3]$. It follows that $\varphi(x)$ is non-decreasing on $[0, 1]$, and hence
\begin{align*}
	{\varphi(x) \le \varphi(1) = (t - 1)(1 - (n\mathbf{r})^2)^2 \le 0, \quad n\mathbf{r} \in [0, 1/3].}
\end{align*}
This completes the proof.
\end{proof}
\section{\bf Concluding remarks} In this paper, we have successfully established several sharp multidimensional analogues of refined Bohr-type inequalities within the unit polydisc $\mathbb{D}^n$. By extending the classical Bohr phenomenon to incorporate the Bohr-Rogosinski radius and Schwarz functions in $\mathcal{B}_{n,m}$, we have demonstrated that the ``$1/3$-phenomenon'' maintains its structural integrity in higher dimensions when appropriately scaled.\vspace{1.2mm}

Our utilization of the radial (Euler) derivative operator $Df(z)$ proved to be the natural vehicle for generalizing univariate derivative growth estimates. Moreover, the resolution of the area-based Bohr inequality in the multivariate setting confirms that the geometric interpretation of the Bohr radius remains consistent across different measure-theoretic perspectives. The sharpness of all obtained results, verified through carefully constructed extremal functions, provides a definitive boundary for these inequalities on $\mathbb{D}^n$.\vspace{2mm}

The results presented herein open several promising avenues for further research:
\begin{prob}(\textbf{Transition to the unit ball:})
	Let $\mathbb{B}^n = \{z \in \mathbb{C}^n : \sum_{j=1}^n |z_j|^2 < 1\}$ be the open unit ball in $\mathbb{C}^n$. Given the sharp radii $R_{m,n,N}$ established for the polydisc $\mathbb{D}^n$ in Theorem \ref{Th-2.2}, a significant open problem is to determine the sharp Bohr-Rogosinski radius $R(\mathbb{B}^n)$ such that for any $f \in \mathcal{H}(\mathbb{B}^n, \mathbb{D})$, the inequality$$|f(\omega(z))| + \sum_{i=1}^{\infty} \sum_{|\alpha|=iN} |a_{\alpha}| r^{\alpha} \leq 1$$holds for $\|z\|_2 \leq R(\mathbb{B}^n)$. Unlike the polydisc, where $|z|^\alpha \leq r^{|\alpha|}$ is applied coordinate-wise, the geometry of $\mathbb{B}^n$ necessitates an investigation into the coupling of variables and the potential dependence of the radius on the norm-induced geometry:$$R(\mathbb{B}^n) = \sup \left\{ r \in (0, 1) : M_f(r) \leq 1, \forall f \in \mathcal{H}(\mathbb{B}^n, \mathbb{D}) \right\}$$where $M_f(r)$ represents the associated Bohr-type sum restricted to the Euclidean ball.
\end{prob}
\begin{prob} (\textbf{General Reinhardt domains:})
	Let $\Omega \subset \mathbb{C}^n$ be a bounded complete Reinhardt domain with Minkowski functional $h_{\Omega}(z) = \inf \{t > 0 : z/t \in \Omega\}$. Given that the polydisc $\mathbb{D}^n$ is a special case where $h_{\mathbb{D}^n}(z) = \max_{1 \le j \le n} |z_j|$, a natural extension of our main results involves determining the sharp Bohr-Rogosinski radius $R(\Omega)$ such that for any $f(z) = \sum_{\alpha} a_{\alpha} z^{\alpha} \in \mathcal{H}(\Omega, \mathbb{D})$:$$|f(\omega(z))| + \sum_{i=1}^{\infty} \sum_{|\alpha|=iN} |a_{\alpha}| \rho_{\alpha} \leq 1 \quad \text{for } h_{\Omega}(z) \leq R(\Omega)$$where $\rho_{\alpha} = \sup_{z \in \Omega} |z^{\alpha}|$. This investigation focuses on how the symmetry and geometric logarithmic convexity of $\partial \Omega$ influence the sharp constants, specifically addressing the interaction between the multi-index coefficients $a_{\alpha}$ and the domain's specific weight distributions
	\begin{align*}
		R(\Omega) := \sup \left\{ r \in (0, 1) : \sup_{f \in \mathcal{H}(\Omega, \mathbb{D})} \left( |f(\omega(z))| + \sum_{\alpha \neq 0} |a_{\alpha} z^{\alpha}| \right) \leq 1,\;z \in r\Omega \right\}.
	\end{align*}
\end{prob}
\begin{prob}(\textbf{Operator-theoretic generalizations and uniform algebras:})
	Let $H^\infty(\mathbb{D}^n)$ denote the Banach algebra of bounded holomorphic functions on the polydisc. A significant extension involves investigating Bohr-type inequalities within the framework of Composition Operators $C_\varphi: f \mapsto f \circ \varphi$. Specifically, one may seek the sharp radius $R_\varphi$ such that for a given self-map $\varphi: \mathbb{D}^n \to \mathbb{D}^n$, the following operator-norm inequality holds:$$\|C_\varphi(f)\|_{\mathcal{B}} := |(f \circ \varphi)(0)| + \sum_{\alpha \neq 0} |a_\alpha| |z^\alpha| \leq \|f\|_\infty \quad \text{for } \|z\|_\infty \leq R_\varphi$$Furthermore, bridging geometric function theory with functional analysis, let $\mathcal{A}$ be a Uniform Algebra on a compact Hausdorff space $X$. We propose the study of the Bohr phenomenon for elements $h \in \mathcal{A}$ by defining a generalized Bohr radius $R(\mathcal{A})$ such that:$$\sup \left\{ \sum_{n=0}^\infty \|h_n\| r^n : h = \sum h_n, \|h\| \leq 1 \right\} \leq 1 \quad \text{for } r \leq R(\mathcal{A})$$where the decomposition $h = \sum h_n$ corresponds to the abstract power series expansion relative to the algebraic structure of $\mathcal{A}$. This would generalize the results from $\mathbb{C}^n$ to a non-commutative or abstract functional setting, as suggested by the properties of Dirichlet series and von Neumann?s inequality.
\end{prob}
\begin{prob} ({\bf Higher-Order Radial Derivatives and Taylor Coefficient Stability}:)
	While this study establishes sharp radii for the first-order Euler operator $Df(z) = \sum_{k=1}^{n} z_k \frac{\partial f}{\partial z_k}$, a significant extension involves the higher-order radial (Euler) operators $D^m f(z)$ defined iteratively as $D^m f(z) = D(D^{m-1} f(z))$ for $m \in \mathbb{N}$. Specifically, we propose the determination of the sharp Bohr-type radius $R_{n,m,\lambda}$ such that for any $f(z) = \sum_{\alpha} a_{\alpha} z^{\alpha} \in \mathcal{H}(\mathbb{D}^n, \mathbb{D})$, the inequality involving the $m$-th order derivative holds:$$|f(z)| + |D^m f(z)| + \lambda \sum_{k=m+1}^{\infty} \sum_{|\alpha|=k} |a_{\alpha}| r^{\alpha} \leq 1 \quad \text{for } \|z\|_{\infty} = r \leq R_{n,m,\lambda}$$Given that the action of the operator on power series is $D^m f(z) = \sum_{\alpha} |\alpha|^m a_{\alpha} z^{\alpha}$, the investigation aims to characterize the growth of the Taylor coefficients $a_{\alpha}$ when weighted by the polynomial factor $|\alpha|^m$. The objective is to find the positive root of the associated majorant series equation:$$A(r) + \sum_{k=1}^{\infty} k^m |a_k| r^k + \lambda \sum_{k=m+1}^{\infty} |a_k| r^k = 1$$extended to the multidimensional setting $z \in \mathbb{C}^n$, which would provide deeper insights into the structural stability of the Bohr phenomenon under repeated differentiation.
\end{prob}

\vspace{5mm}

\noindent\textbf{Conflict of interest:} The authors declare that there is no conflict  of interest regarding the publication of this paper.\vspace{1.2mm}

\noindent {\bf Funding:} Not Applicable.\vspace{1.2mm}

\noindent\textbf{Data availability statement:}  Data sharing not applicable to this article as no datasets were generated or analysed during the current study.\vspace{1.2mm}

\noindent {\bf Authors' contributions:} All the authors have equal contributions in preparation of the manuscript.

\end{document}